\documentclass[10pt,twoside]{article}
\usepackage[latin1]{inputenc}
\usepackage{amsmath}
\usepackage{graphicx}
\usepackage{yhmath}
\usepackage[table]{xcolor}
\usepackage{mathrsfs} 
\usepackage{amssymb}
\usepackage{url}
\usepackage{makecell}
\usepackage{arydshln}
\usepackage{multirow}
\usepackage{arydshln}
\usepackage{mathtools}
\usepackage{stmaryrd}  

\input epsf
\def\limiten{\renewcommand{\arraystretch}{0.5}
\begin{array}[t]{c}\stackrel{}{\longrightarrow} \\
{\scriptstyle n\rightarrow
\infty}\end{array}\renewcommand{\arraystretch}{1}}

\def\limitepsn{\renewcommand{\arraystretch}{0.5}
\begin{array}[t]{c}\stackrel{a.s.}{\longrightarrow} \\
{\scriptstyle n \rightarrow
\infty}\end{array}\renewcommand{\arraystretch}{1}}

\def\limiteloin{\renewcommand{\arraystretch}{0.5}
\begin{array}[t]{c}\stackrel{{\cal D}}{\longrightarrow} \\
{\scriptstyle n\rightarrow
\infty}\end{array}\renewcommand{\arraystretch}{1}}

\def\limiteproban{\renewcommand{\arraystretch}{0.5}
\begin{array}[t]{c}\stackrel{{\cal P}}{\longrightarrow} \\
{\scriptstyle n\rightarrow
\infty}\end{array}\renewcommand{\arraystretch}{1}}

\numberwithin{equation}{section}

\newtheorem{thm}{Theorem}[section]

\newtheorem{Corol}[thm]{Corollary}

\newtheorem{lem}[thm]{Lemma}

\newtheorem{prop}[thm]{Proposition}

\newcommand{\E}{\ensuremath{\mathbb{E}}}
\newcommand{\R}{\ensuremath{\mathbb{R}}}
\newcommand{\Z}{\ensuremath{\mathbb{Z}}}

\newcommand{\N}{\ensuremath{\mathbb{N}}}

\definecolor{grisclair}{gray}{0.9}

\renewcommand{\arraystretch}{.8}

\begin{document}
\title{\bf Inference and model selection in general causal time series with  exogenous covariates}
 \maketitle \vspace{-1.0cm}
 \begin{center}
   Mamadou Lamine DIOP \footnote{Supported by
   the MME-DII center of excellence (ANR-11-LABEX-0023-01) 
   } 
   and 
    William KENGNE \footnote{Developed within the ANR BREAKRISK: ANR-17-CE26-0001-01 and the  CY Initiative of Excellence (grant "Investissements d'Avenir" ANR-16-IDEX-0008), Project "EcoDep" PSI-AAP2020-0000000013} 
 \end{center}

  \begin{center}
  { \it 
 THEMA, CY Cergy Paris Université, 33 Boulevard du Port, 95011 Cergy-Pontoise Cedex, France.\\
  E-mail: mamadou-lamine.diop@u-cergy.fr ; william.kengne@u-cergy.fr  \\
  }
\end{center}

 \pagestyle{myheadings}
 \markboth{Inference and model selection in general causal time series with  exogenous covariates}{Diop and Kengne}

~~\\
\textbf{Abstract}:
 In this paper, we study a general class of causal processes with exogenous covariates, including  many classical processes such as the ARMA-GARCH, APARCH, ARMAX, GARCH-X and APARCH-X processes.
 Under some Lipschitz-type conditions, the existence of a $\tau$-weakly dependent strictly stationary and ergodic solution is established. 
 We provide conditions for the strong consistency and derive the asymptotic distribution of the quasi-maximum likelihood estimator (QMLE), both when the true parameter is an interior point of the parameters space and when it belongs to the boundary. 
 A significance Wald-type test of parameter is developed. This test is quite extensive and includes the test of nullity of the parameter's components, which in particular, allows us to assess the relevance of the exogenous covariates.
 Relying on the QMLE of the model, we also propose a penalized criterion to address the problem of the model selection for this class.  The weak and the strong consistency of the procedure are established.
 Finally, Monte Carlo simulations are conducted to numerically illustrate the main results.
 \\ 


\section{Introduction}
Autoregressive time series with exogenous covariates provide effective ways to take into account some available extra information in the models. The well known example that has been widely studied is the ARMAX model, see Hannan (1976), Hannan and Deistler (2012). 
The GARCH-type models with exogenous covariates have recently attracted much attention in the literature, see for instance Han and Kristensen (2014) for GARCH-X, Francq and Thieu (2019) for  APARCH-X. Guo {\it et al.} (2014) considered the factor double autoregressive model, whose ARX and ARCH-X are particular cases. We consider a large class of causal time series models, whose ARMAX and GARCH-X type models are specific examples.

\medskip

 Let $X_t=(X_{1,t},X_{2,t},\ldots,X_{d_{x}, t}) \in \R^{d_{x}}$ be a vector of covariates, with $d_{x} \in \N$.  Consider the class of affine causal models with exogenous covariates,
 
 \medskip
 
 \noindent \textbf{Class} $\mathcal {AC}$-$X(M_\theta,f_\theta):$ A process $\{Y_{t},\,t\in \Z \}$ belongs to $\mathcal {AC}$-$X(M_\theta,f_\theta)$ if it satisfies:
   \begin{equation}\label{Model} 
     Y_t =M_\theta(Y_{t-1}, Y_{t-2}, \ldots; X_{t-1},X_{t-2},\ldots)\xi_t + 
      f_\theta(Y_{t-1}, Y_{t-2}, \ldots; X_{t-1},X_{t-2},\ldots),
   \end{equation}
  where 
   $M_\theta,~f_\theta :  \R^{\N} \times (\R^{d_x})^{\N} \rightarrow \R$ are two measurable functions 
   and assumed to be known up to the parameter $\theta$, which belongs in a compact subset $ \Theta\subset \R^d$ ($d \in \N$); 
   and
   $(\xi_t)_{t \in \Z}$ is a sequence of zero-mean independent, identically distributed (\textit{i.i.d}) random variable satisfying $\E(\xi^r_0) < \infty$ for some $r \geq 2$ and $\E(\xi^2_0) =1$. 
   Remark that, if $X_t \equiv C$ for some constant $C$ (absence of covariates), then (\ref{Model}) reduces to the classical affine causal models that has already been considered in the literature (see, for instance, Bardet and Wintenberger (2009), Bardet {\it et al.} (2012), Bardet {\it et al.} (2020)). One can see that, the  ARMAX, GARCH-X, APARCH-X models belong to the class $\mathcal {AC}$-$X(M_\theta,f_\theta)$. 
   
    \medskip
 
 There exist several important contributions devoted to autoregressive models with covariates; we refer to Hannan and Deistler (2012), Han and Kristensen (2014),  Sucarrat {\it et al.} (2016),  Francq and Sucarrat (2017), Pedersen and Rahbek (2018), Francq and Thieu (2019), Gr{\o}nneberg and Holcblat (2019),   Zambom and Gel (2020) and the references therein for some developments on ARMAX and conditional volatility type models with exogenous covariates. The class $\mathcal {AC}$-$X(M_\theta,f_\theta)$ is more general than the models considered in the aforementioned works, as well as the factor double autoregressive model proposed by Guo {\it et al.} (2014) which is a particular case of the model (\ref{Model}). Note as well that, the class $\mathcal {AC}$-$X(M_\theta,f_\theta)$ provides a more general way to take into account covariates in the model, and one can see that the linear covariates regressors considered by Francq and Thieu (2019) and many other works is a specific case.  
 Compared to Bardet and Wintenberger (2009), besides taking into account covariates in the model (\ref{Model}), we address the inference when the true parameter belongs to the boundary of the parameter set $\Theta$ and the model selection question. 
 
  \medskip
  
  In this new contribution, we consider the class of model (\ref{Model}) and  address the following issues.
  \begin{enumerate}
  \item[(i)] \textbf{Existence of a stationary solution}. We provide sufficient conditions that ensure the existence of a $\tau$-weakly dependent stationary and ergodic solution  $Z_t=(Y_t,X_t)$ of (\ref{Model}). At a first glance, one might think that these conditions are the same as those obtained by Bardet and Wintenberger (2009), but in our case, the existence of the covariates must be taken into account.
  \item[(ii)] \textbf{Inference for the class $\mathcal {AC}$-$X(M_{\theta^*},f_{\theta^*})$}. An inference based on the quasi likelihood of the model is carried out. The consistency of the quasi-maximum likelihood estimator (QMLE) is established and we derived the asymptotic distribution of this estimator (even when $\theta^*$  belongs to the boundary of $\Theta$). 
  \item[(iii)] \textbf{Significance test of parameter}. A Wald-type significance test of parameter of the model (\ref{Model}) is conducted. The proposed test is quite extensive and includes the test of nullity of the parameter's components. An asymptotic study is carried out, which shows in particular that, when the true parameter belongs to the boundary of $\Theta$, the asymptotic distribution of the test statistic under the null hypothesis is quite different from the classical chi-square distribution.
  \item[(iv)] \textbf{Model selection}. A penalized criterion based on the quasi likelihood of the model is proposed for model selection in the class $\mathcal {AC}$-$X(M_{\theta^*},f_{\theta^*})$. We provides conditions that ensure the weak and the strong consistency of the proposed procedure. These conditions shows in particular that, the Hannan-Quinn information Criterion (HQC) with a regularization parameter $\kappa_n = c \log \log n$ (see (\ref{Cont_pen})) is strongly consistent for sufficiently large $c$.
  \end{enumerate}
%
%
The article is organized as follows. 
In Section 2, firstly, we provide conditions for stability properties. Secondly, we give the definition of the QMLE and study its asymptotic properties; a significance test of parameter with an asymptotic study is also addressed.
Section 3 focuses on the model selection and the consistency of the proposed procedure. Some classical examples of processes belonging to the class $\mathcal {AC}$-$X(M_{\theta^*},f_{\theta^*})$ are detailed in Section 4.
Section 5 gives some empirical results, whereas Section 6 is devoted to a summary and conclusion.
Section 7 contains the proofs of the main results.

\section{Assumptions, inference and  test of the parameters}
\subsection{Assumptions}
 
    \noindent
 Throughout the sequel, the following norms will be used:
 {\em
\begin{itemize}
 \item $ \|x \| \coloneqq \sqrt{\sum_{i=1}^{p} x_i^2 } $, for any $x \in \mathbb{R}^{p}$, $p \in \N$;
%
\item $ \|V\| \coloneqq \sqrt{\sum_{i=1}^{p} \sum_{j=1}^{q} v_{i,j}^2 } $, for any matrix $V \in M_{p,q}(\R)$,  where $M_{p,q}(\R)$ denotes the set of matrices of dimension $p\times q$ with coefficients in $\R$, for $p,q \in \N$;
\item  $\left\|g\right\|_{\mathcal K} \coloneqq \sup_{\theta \in \mathcal K}\left(\left\|g(\theta)\right\|\right)$ for any compact set $\mathcal K \subseteq \R^d$ and function $g:\mathcal K \longrightarrow   M_{p,q}(\R)$;
 %
%
\item $\left\|Y\right\|_r \coloneqq \E\left(\left\|Y\right\|^r\right)^{1/r}$, if $Y$ is a random vector with finite $r-$order moments, for $r >0$.
\end{itemize}
}

 \medskip
 
  \noindent
  We will denote by $0$ the null vector of any vector space. Let $\Psi_\theta$ be the generic symbol for any of the functions $f_\theta$ or $M_\theta$. 
We set the following classical Lipschitz-type conditions for any compact set $\mathcal K \subseteq \Theta$.

   \medskip
   
    \noindent 
    \textbf{Assumption} \textbf{A}$_i (\Psi_\theta,\mathcal K)$ ($i=0,1,2$):
    For any $(y,x) \in \R^{\N}\times (\R^{d_x})^{\N}$, the function $\theta \mapsto \Psi_\theta(y)$ is $i$ times continuously differentiable on $\mathcal K$  with $ \big\| \frac{\partial^i \Psi_\theta(0)}{\partial \theta^i}\big\|_{\mathcal K}<\infty $; 
    and
      there exists two sequences of non-negative real numbers $(\alpha^{(i)}_{k,Y}(\Psi_\theta,\mathcal K))_{k\geq 1} $ and $(\alpha^{(i)}_{k,X}(\Psi_\theta,\mathcal K))_{k\geq 1}$ satisfying:
     $ \sum\limits_{k=1}^{\infty} \alpha^{(i)}_{k,Y}(\Psi_\theta,\mathcal K) <\infty$, $ \sum\limits_{k=1}^{\infty} \alpha^{(i)}_{k,X}(\Psi_\theta,\mathcal K) <\infty$  for $i=0, 1, 2$;
   such that for any  $(y,x), (y',x') \in \R^{\N} \times (\R^{d_x})^{\N}$,
  \[  \Big \| \frac{\partial^i \Psi_\theta(y,x)}{ \partial \theta^i}-\frac{\partial^i \Psi_\theta(y',x')}{\partial\theta^i} \Big \|_{\mathcal K}
  \leq  \sum\limits_{k=1}^{\infty}\alpha^{(i)}_{k,Y}(\Psi_\theta,\mathcal K) |y_k-y'_k|+\sum\limits_{k=1}^{\infty}\alpha^{(i)}_{k,X}(\Psi_\theta,\mathcal K)\|x_k-x'_k\|), 
  \]
where $\| \cdot\|$ denotes any vector, matrix norm.
\medskip
   
    \noindent
The following assumption is considered on the function $H_\theta=M^2_\theta$ in the cases of ARCH-X type process.
\medskip
   
    \noindent
 \textbf{Assumption} \textbf{A}$_i (H_\theta,\mathcal K)$ ($i=0,1,2$): Assume that $f_\theta=0$. 
   There exists two sequences of non-negative real numbers $(\alpha^{(i)}_{k,Y}(H_\theta,\mathcal K))_{k\geq 1} $ and $(\alpha^{(i)}_{k,X}(H_\theta,\mathcal K))_{k\geq 1} $ satisfying:
     $ \sum\limits_{k=1}^{\infty} \alpha^{(i)}_{k,Y}(H_\theta,\mathcal K) <\infty $, $ \sum\limits_{k=1}^{\infty} \alpha^{(i)}_{k,X}(H_\theta,\mathcal K) <\infty $  for $i=0, 1, 2$;
   such that for any  $(y,x), (y',x') \in \R^{\infty} \times (\R^{d_x})^{\infty}$,
  \[  \Big \| \frac{\partial^i H_\theta(y,x)}{ \partial \theta^i}-\frac{\partial^i H_\theta(y',x')}{\partial\theta^i} \Big \|_{\mathcal K}
  \leq  \sum\limits_{k=1}^{\infty}\alpha^{(i)}_{k,Y}(H_\theta,\mathcal K) |y^2_k-{y'}^2_k|+\sum\limits_{k=1}^{\infty}\alpha^{(i)}_{k,X}(H_\theta,\mathcal K) \|x_k-x'_k\| . 
  \]
  
\medskip

  \noindent
In the whole paper, we impose an autoregressive-type structure on the covariates:
 \begin{equation}\label{Process_X}
 X_t=g(X_{t-1},X_{t-2},\ldots;\eta_t),
\end{equation}
 where $(\eta_t)_{t \in \Z}$ is a sequence of zero-mean random variables such as $(\eta_t, \xi_t)_{t \in \Z}$  is {\it  i.i.d}
  and $g$ is a function with values in $\R^{d_x}$ satisfying
\begin{equation} \label{exp_Lip_cov}  
   \E\left[\left\|g(0, \eta_0) \right\|^r\right]<\infty
   ~ \text{ and }  ~ 
   \left\|g(x; \eta_0)-g(x'; \eta_0) \right\|_r \leq \sum\limits_{k=1}^{\infty} \alpha_k(g) \left\| x_k-x'_k \right\|
   ~ \text{ for all } x, x'  \in (\R^{d_{x}})^{\N},
  \end{equation}
   for some $r \geq 1$ and non-negative sequence $(\alpha_k(g))_{k \geq 1}$ such that $\sum\limits_{k=1}^{\infty} \alpha_k(g) <1 $.
   
   \medskip
   
    \noindent 
    For $r\geq 1$, when (\ref{exp_Lip_cov}) holds,  we define the set 
\begin{multline*}\label{Set_Theta(r)}
\Theta(r) = \Big\{
 \theta \in \R^d \, \big / \, \textbf{A}_0 (f_\theta,\{\theta\}) \   \text{and}\ \textbf{A}_0 (M_\theta,\{\theta\})   \    \text{hold with} \\ 
 \hspace{7.01cm} \sum\limits_{k=1}^{\infty} \max \left\{\alpha_k(g),\, \alpha^{(0)}_{k,Y}(f_\theta,\{\theta\}) + \|\xi_0\|_r \alpha^{(0)}_{k,Y}(M_\theta,\{\theta\}\right\} 
 <1
\Big\}\\
\bigcup 
\Big\{
 \theta  \in \R^d \ \big / \ f_\theta=0 \text{ and } \textbf{A}_0 (H_\theta,\{\theta\})  \text{ holds with } 
 \|\xi_0\|^2_r  \sum\limits_{k=1}^{\infty} \max \left\{\alpha_k(g),\,\alpha^{(0)}_{k,Y}(H_\theta,\{\theta\}) \right\}
 <1
\Big\}.
\end{multline*}

\medskip

  \noindent
In the sequel, we make the convention that if \textbf{A}$_i(M_\theta,\Theta)$ holds then $\alpha^{(i)}_{k,Y}(H_\theta,\Theta)=\alpha^{(i)}_{k,X}(H_\theta,\Theta)=0$ for all $k\in\N$ and if 
\textbf{A}$_i(H_\theta,\Theta)$ holds then $\alpha^{(i)}_{k,Y}(M_\theta,\Theta)=\alpha^{(i)}_{k,X}(M_\theta,\Theta)=0$ for all $k\in\N$.

   \medskip
   
    \noindent 
    The condition (\ref{exp_Lip_cov}) allows to assure the stability of the process $X_t$. 
    Together the aforementioned assumptions assure the existence of a stationary and weakly dependent solution  of order $r$ to the model (\ref{Model}), as shown in the following proposition.
    \begin{prop}\label{prop1}
    Assume that  \textbf{A}$_0(f_\theta,\Theta)$, \textbf{A}$_0(M_\theta,\Theta)$ (or \textbf{A}$_0(H_\theta,\Theta)$) and (\ref{exp_Lip_cov}) hold. If $\theta^* \in \Theta \cap \Theta(r)$ with $r\geq 1$, then there exists a $\tau$-weakly dependent stationary, ergodic and non anticipative solution $(Z_t)_{t\in \Z}$ $Z_t=(Y_t,X_t)$, to (\ref{Model}), satisfying $ \E[{\| Z_0\|}^r] <\infty$. 
    \end{prop}
  
\subsection{Inference and significance test of parameter}
In this paragraph, we describe the use of the Gaussian quasi-maximum likelihood to obtain an estimator of the parameters of the model (\ref{Model}). The  main asymptotic properties of this estimator are also established.    
Assume that the observations $(Y_1,X_1),\ldots,(Y_n,X_n)$ are generated from (\ref{Model}) and (\ref{Process_X}) according to the true parameter $\theta^* \in \Theta$ which is unknown. For all $t \in \Z$, denote by $\mathcal{F}_{t}=\sigma((Y_s,X_s),\, s \leq t)$ the $\sigma$-field generated by the whole past at time $t$. 
%
%
%
 The mean and the variance of $Y_t | \mathcal{F}_{t-1}$ and $f_{\theta^*}(Y_{t-1},\ldots;X_{t-1},\ldots)$ and variance $M^2_{\theta^*}(Y_{t-1},\ldots;X_{t-1},\ldots)$ respectively. 
For any $\theta \in \Theta$, the conditional Gaussian quasi-log-likelihood is given by (up to an additional constant)

\[  L_n(\theta) :=  -\frac{1}{2}\sum_{t=1}^{n} q_t(\theta) ~~~ \text{ with }~~~
        q_t(\theta) = \frac{ (Y_t-f^t_\theta)^2}{H^t_\theta}+\log H^t_\theta,
  \]
  where $f^t_\theta := f_\theta(Y_{t-1},Y_{t-2}\ldots;X_{t-1},X_{t-2},\ldots)$,  
  $M^t_\theta := M_\theta(Y_{t-1},Y_{t-2}\ldots;X_{t-1},X_{t-2},\ldots)$ and $H^t_\theta := (M^t_\theta)^2$.\\
Since $(Y_0,X_0),(Y_{-1},X_{-1}),\ldots$ are not observed, $L_n(\theta)$ is approximated by 
\[ 
 \widehat  L_n(\theta)=  -\frac{1}{2}\sum_{t=1}^{n} \widehat q_t(\theta) ~~~ \text{ with }~~~
        \widehat q_t(\theta) = \frac{ (Y_t-\widehat f^t_\theta)^2}{\widehat H^t_\theta}+\log \widehat H^t_\theta,
 \]       
where $\widehat f^t_\theta := f_\theta(Y_{t-1},\ldots,Y_{1},0;X_{t-1},\ldots,X_{1},0)$, 
$\widehat M^t_\theta := M_\theta(Y_{t-1},\ldots,Y_{1},0;X_{t-1},\ldots,X_{1},0)$ and $\widehat H^t_\theta := (\widehat M^t_\theta)^2$.
 Thus, the QMLE of $\theta^*$ is defined by
   \[
 \widehat{\theta}_{n}= \underset{\theta\in \Theta}{\text{argmax}} \big( \widehat L_n(\theta) \big).
 \]
 %
          
 \noindent
 We set the following regularity conditions to assure the  identifiability of the model and to derive the asymptotic behavior of the QMLE. 
 \begin{enumerate}
  \item [(\textbf{A0}):] for all  $\theta \in \Theta$ and some $t \in \Z$, 
 $ \big( f^t_{\theta^*}= f^t_{\theta}  \ \text{and} \ H^t_{\theta^*}= H^t_{\theta} \ \ a.s. \big) \Rightarrow ~ \theta= \theta^*$;

  \item [(\textbf{A1}):] $\exists  \underline{h}>0$ such that $\displaystyle \inf_{ \theta \in \Theta} H_\theta(y,x)  \geq \underline{h}$,  for all $(y,x) \in \R^{\N}\times (\R^{d_x})^{\N}$;
  
  
    \item [(\textbf{A2}):] for all  $\theta \in \Theta$, $c \in \R^d$, $\Big($ $c' \frac{\partial}{\partial \theta}f^0_{\theta^*}=0$ or $c' \frac{\partial}{\partial    \theta}H^0_{\theta^*}=0$ $\Big) ~a.s.$   $\Longrightarrow~ c=0$,  where $'$ denotes the transpose.
\end{enumerate}

  \medskip  

\noindent
Assumption (\textbf{A0}) is an identifiability condition and it will be discussed in detail for each of the examples of processes studied in the paper. From (\textbf{A1}), the quasi likelihood is well defined, whereas (\textbf{A2}), which is classical (see for instance Bardet and Wintenberger (2009)) allows to derive the asymptotic distribution of the QMLE. The following theorem addresses the strong consistency of the QMLE.
\begin{thm}\label{th1}
Assume that (\textbf{A0}), (\textbf{A1}), \textbf{A}$_0(f_\theta,\Theta)$, \textbf{A}$_0(M_\theta,\Theta)$ 
and (\ref{exp_Lip_cov}) (with $r \geq 2$) hold with 
\begin{equation}\label{cond_th1}
 \alpha^{(0)}_{k,Y}(f_\theta,\Theta)+\alpha^{(0)}_{k,X}(f_\theta,\Theta)+\alpha^{(0)}_{k,Y}(M_\theta,\Theta)+\alpha^{(0)}_{k,X}(M_\theta,\Theta)
 +\alpha^{(0)}_{k,Y}(H_\theta,\Theta)+\alpha^{(0)}_{k,X}(H_\theta,\Theta) = O(k^{-\gamma}),
\end{equation}
 for some $\gamma >3/2$. \\
If $\theta^* \in \Theta \cap \Theta(r) $ with $r \geq 2$, then  
\[
\widehat{\theta}_{n} \limitepsn \theta^*.
\]
\end{thm}

\medskip
 
 \medskip
 \noindent
%
To derive the asymptotic distribution of the QMLE, it is necessary to take into account the constraints in the parameter space $\Theta$ corresponding to the model. 
 For example, in some processes belonging to (\ref{Model}), such as the ARCH-X models (see below), the  components of $\theta^*$ are constrained to be positive or equal to zero. 
 In order to propose a parsimonious representation, it is often required to test whether or not the exogenous covariates are relevant.
 For example, in an ARCH-X$(1)$ model defined by 
 $Y_t=  \xi_{t} \sigma_t 
  \text{ with } 
 \sigma^2_t=\alpha^*_0 + \alpha^*_1 Y^2_{t-k}+ {\gamma^*}'X_{t-1}$, the true parameter vector is $\theta^* = (\alpha^*_0, \alpha^*_1,\gamma^*) \in \Theta \subset ]0,\infty[ \times {[0,\infty[}^{dx+1}$. The significance test of the covariate $X_t$ consists to verify  the nullity of the parameter $\gamma^*$; that is, 
 if the true parameter vector can be of the form $\theta_0=(\alpha_0, \alpha_1,0)$ which is not an interior point of $\Theta$.
In this situation, 
 it is impossible to apply the asymptotic normality results based on the classical assumption of "interior point" to derive the asymptotic behavior of the test statistic used. 
To take into account such a scenario in the general class (\ref{Model}), we will consider that the component $i$ of  $\theta^*$ is constrained if the $i$-th section of $\Theta$ is of the form   $[\underline{\theta}_i, \overline{\theta}_i]$ with $\underline{\theta}_i< \overline{\theta}_i$. 
Assume that the   $d_2$ (with $d_2 \in \{ 0,\ldots,d\}$) last components of $\theta^*$ are constrained, and let $d_1=d-d_2$.  
Therefore,  if $d_2 \geq 1$ and  $\theta^*_i \in \{ \underline{\theta}_i, \overline{\theta}_i\}$ with $i >d_1$, then $\theta^*$ is not an interior point of $\Theta$. For instance, in a scenario where $\theta^*_i= \underline{\theta}_i$,
 with the QMLE $\widehat{\theta}_{n}=(\widehat \theta_{1,n}, \ldots,\widehat \theta_{d,n})$, it holds that $\sqrt{n}(\widehat \theta_{i,n}-\theta^*_{i})\in [0,\infty)$ which cannot tend to a Gaussian distribution  with mean $0$.
By convention, it is assumed that $\theta^* \in \overset{\circ}{\Theta}$ if  $d_2=0$. 
 When $d_2 \geq 1$ and the set $\Theta$ is assumed to be large enough, then the following relation holds:
\begin{equation}\label{eq0_normal}
 \bigcup_{n \geq 1}\left\{\sqrt{n}\left(\theta-\theta^*\right),\ \theta \in \Theta \right\}= \mathcal C ~~\text{ with } ~~ \mathcal C= \prod_{i=1}^{d}\mathcal C_i, 
\end{equation}
where
$\mathcal C_i = [0, \infty[$  when $i>d_1$ and $\theta^*_i = \underline{\theta}_i$, $\mathcal C_i = ]\infty, 0]$  when $i>d_1$ and $\theta^*_i = \overline{\theta}_i$,  and  $\mathcal C_i=\R$ otherwise. 
The set $\mathcal C$ is a convex cone which is equal to $\R^d$ if $\theta^* \in \overset{\circ}{\Theta}$. 

 \medskip
 
\noindent
Let us define the following matrices
\begin{equation}\label{def_matrix_F_G}
 F= \E\Big[\frac{ \partial^2  q_{0}(\theta^*)}{\partial \theta \partial \theta'}\Big]
 ~~~ \text{ and } ~~~
  G= \E\Big[ \frac{ \partial  q_{0}(\theta^*)}  {\partial \theta} \frac{ \partial  q_{0}(\theta^*)}{\partial \theta'}  \Big].
\end{equation}
Under the  assumptions \textbf{A}$_i (f_\theta,\Theta)$, \textbf{A}$_i (M_\theta, \Theta)$ (with $i=0,1,2$), one can  show the existence of $F$ and $G$. In addition, in view to (\textbf{A2}), the same arguments as in Bardet and Wintenberger (2009) allow to establish that matrix $F$ is positive definite.    
 Consider then the $F$-scalar product $\langle x,y\rangle_F =x'Fy$ and the norm ${\|x \|}^2_F= x'Fx$  for $x,y \in \R^d$. 
Let us define the $F$-orthogonal projection of a vector $Z \in \R^d$ on the cone $\mathcal C$ as follows:
\begin{equation*}
Z^{\mathcal C} = \text{arg} \inf_{C \in \mathcal C} \left\| C-Z\right\|_F.
\end{equation*}
This definition  is equivalent to
\begin{equation}\label{def_project_orth}
Z^{\mathcal C} \in  \mathcal C ~~\text{ with } ~~\left\langle Z-Z^{\mathcal C}, C-Z^{\mathcal C}\right\rangle \leq 0,~ \ \forall C \in \mathcal C.
\end{equation}
Note that, when $\theta^* \in \overset{\circ}{\Theta}$, we have $Z^{\mathcal C} =Z$.  
%
Combining all the regularity conditions and definitions given above, we obtain the following main result.

\begin{thm}\label{th2}
Assume that (\textbf{A0})-(\textbf{A2}),  (\textbf{A}$_i(f_\theta,\Theta)$), (\textbf{A}$_i(M_\theta,\Theta)$) (for $i=0,1,2$) and (\ref{exp_Lip_cov}) (with $r \geq 4$) hold with
\begin{equation}\label{cond_th2} \alpha^{(i)}_{k,Y}(f_\theta,\Theta)+\alpha^{(i)}_{k,X}(f_\theta,\Theta)+\alpha^{(i)}_{k_\theta,Y}(M_\theta,\Theta)+\alpha^{(i)}_{k,X}(M_\theta,\Theta)
   +\alpha^{(i)}_{k_\theta,Y}(H_\theta,\Theta)+\alpha^{(i)}_{k,X}(H_\theta,\Theta) = O(k^{-\gamma}) 
  ,
\end{equation}
 for  $i=0,1,2$  and some $\gamma >3/2$.
\begin{itemize}
	\item If $\theta^* \in \Theta \cap \Theta(r)$ with  $r \geq 4$, then
 \[
 \sqrt{n}\left(\widehat\theta_{n}-\theta^*\right) \limiteloin Z^{\mathcal C} ~~~ with ~~~
  Z\sim \mathcal{N}_d \left(0,\Sigma \right),
 ~~~where~~~ \Sigma :=F^{-1} G F^{-1}.
 \]
 
 \item  If $\theta^* \in \overset{\circ}{\Theta} \cap \Theta(r)$ with  $r \geq 4$, then 
 \[
 \widehat\theta_{n} - \theta^* =O\bigg(\sqrt{\frac{\log \log n}{n}}\bigg)~~ a.s.
 \]
 \end{itemize}
 \end{thm}
 
 \medskip
 
 \medskip
 \noindent
 The matrix $\Sigma$ can be consistently estimated by $\widehat \Sigma_n = F_n(\widehat\theta_{n})^{-1} G_n(\widehat\theta_{n}) F_n(\widehat\theta_{n})^{-1}$, where
 \[
 F_n(\widehat\theta_{n})= \frac{1}{n} \sum_{t=1}^{n}\frac{ \partial^2  q_{t}(\widehat\theta_{n})}{\partial \theta \partial \theta'} 
 ~~\text{ and }~~
G_n(\widehat\theta_{n})= \frac{1}{n} \sum_{t=1}^{n}\frac{ \partial  q_{t}(\widehat\theta_{n})}  {\partial \theta} \frac{ \partial  q_{t}(\widehat\theta_{n})}{\partial \theta'} .
 \]

 \medskip  

\noindent
Now, we are interested to investigate whether or not a given subset of components of $\theta^*$ are equal to some fixed vector. 
 To do so, consider the following hypothesis testing:
 \begin{equation}\label{Test}
 H_0: \, \Gamma \theta^*=\vartheta_0  ~~~~~~ \text { against } ~~~~~~  H_1: \, \Gamma \theta^* \neq \vartheta_0,
\end{equation}
 where $\Gamma $ is a $d_0\times d$ full-rank matrix and $\vartheta_0 $ is a vector of dimension $d_0$. 
 Define the  Wald-type test statistic given by 
 \begin{equation}\label{Stat_Wn}
 W_n= n(\Gamma \widehat\theta_{n}-\vartheta_0)' (\Gamma  \widehat \Sigma_n \Gamma ')^{-1}(\Gamma \widehat\theta_{n}-\vartheta_0).
 \end{equation}
 Under $H_0$, the  asymptotic behavior of $W_n$ is given by the following theorem.
 \begin{thm}\label{th3}
Under $H_0$, assume that the assumptions of Theorem \ref{th2} hold. Then
\[
W_n \limiteloin (\Gamma Z^{\mathcal C})'(\Gamma \Sigma \Gamma')^{-1}\Gamma Z^{\mathcal C} ~~~~\text {with}~~~~ Z \sim \mathcal{N}_d \left(0,\Sigma  \right).
\]
 \end{thm}
By the above theorem, at a nominal level $\alpha \in (0,1)$, the critical region of the test is $(W_{n}>q_\alpha)$, where $q_\alpha$ is the $(1-\alpha)$-quantile of the distribution of $(\Gamma Z^{\mathcal C})'(\Gamma \Sigma \Gamma')^{-1}\Gamma Z^{\mathcal C}$. The critical value $q_\alpha$ can be computed through Monte-Carlo simulations. 
The following corollary follows immediately when $\theta^*$ belongs to the interior of the parameter space. 
\begin{Corol}\label{corol1}
Assume that the conditions of Theorem \ref{th3} hold. If $\theta^* \in \overset{\circ}{\Theta}$, then $W_n$ converges to a chi-square distribution with $d_0$  degrees of freedom.
 \end{Corol}
Under $H_1$, one can easily see that $W_n  \underset{n \rightarrow \infty}{\overset{a.s}{\longrightarrow}}
 +\infty$; which shows that the test is consistent in power. 
In the empirical studies, we will restrict our attention to test the relevance of the exogenous covariates by using the hypothesis (\ref{Test}) with $\vartheta_0=0$ and an appropriate matrix $\Gamma$.

\section{Model selection}\label{Model_selection}

\subsection{Model selection framework} 
Assume that $(Y_1,\ldots,Y_n)$ is a trajectory of the process $Y=\{Y_t,\, t\in \Z \}$ satisfying $\mathcal {AC}$-$X(M_{\theta^*},f_{\theta^*})$ (defined as in (\ref{Model})), where the true parameter $\theta^*$ is unknown. 
Let $\mathcal M$ be a finite collection of models belonging to  $\mathcal {AC}$-$X(M_\theta,f_\theta)$ with $\theta \in \Theta$.  Assume that $\mathcal M$ contains at least the true model $m^*$ corresponding to the parameter $\theta^*$.
Our objective is to develop a procedure that allows to select the "best model" (that we denote by  $\widehat m_n$) among the collection $\mathcal M$ such that it is "close" to $m^*$ for $n$ large enough. 
To this end, we consider the following definitions and notations in the sequel: 
\begin{itemize}
	\item a model $m \in \mathcal M$ is considered as a subset of $\{ 1,\ldots,d\}$ and denote by $|m|$ the dimension of  $m$ (i.e, $|m|=\#(m)$);
	%
		\item for  $m \in \mathcal M$, $\Theta_m=\{  (\theta_i)_{1 \leq i \leq d} \in \Theta \text{ with } \theta_i=0 \text{ if } i \notin m \}$ is  a compact set containing $\theta(m)$, where $\theta(m)$ denotes the parameter vector associated to the model $m$; 
		\item $\mathcal{M}$ is considered as a subset of the power set of $\{ 1,\ldots,d\}$; that is, $\mathcal{M} \subset \mathcal{P}(\{ 1,\ldots,d\})$.
\end{itemize}
For instance, when the observations $Y_1,\ldots,Y_n$ are generated from a ARMAX$(p^*,q^*,s^*)$ model (defined below), 
the collection $\mathcal{M}$ of the competing models could be considered as a family of ARMAX$(p,q,s)$ with $(p,q,s) \in \{0,1,\ldots, p_{max} \} \times \{0,1,\ldots, q_{max} \}\times \{0,1,\ldots, s_{max} \}$, where $p_{max}, q_{max}, s_{max}$ are the fixed upper bounds of the orders satisfying $p_{max} \geq p^*$, $q_{max} \geq q^*$, $s_{max} \geq s^*$.
 The parameter space $\Theta$ is a compact subset of $ \R^{p_{max} + q_{max}+ s_{max}}$, and thus  a model $m$ is a subset of $\{ 1,2,\ldots,p_{max} + q_{max}+ s_{max}\}$.
 
 \subsection{Model selection criterion and asymptotic results}
Note that, under the identifiability assumption (\textbf{A0}), one can show that, for all $m \in \mathcal{M}$, the function  $\theta \mapsto -\E[q_{0}(\theta)]$ has a unique maximum in $\Theta_m$ (see proof of Theorem \ref{th1}). Let us thus define the
"best" parameter associated to the model $m$ as 
\[
\theta^*(m):= \underset{\theta\in \Theta_m}{\text{argmin}} \left(\E[q_{0}(\theta)]\right).
\]  
When $m \supseteq m^*$, we have  $\theta^*(m)=\theta^*(m^*)=\theta^*$; that is, $\theta^*(m)$ will play the role of the true parameter $\theta^*$ in cases of "true" or overfitted model.
For $m\in \mathcal M$, we define the QMLE of $\theta^*(m)$ as 
 \begin{equation}\label{qmle_m}
  \widehat{\theta}(m) :=  \underset{\theta\in \Theta_m}{\text{argmax}} \left(\widehat{L}_n(\theta)\right).
  \end{equation}
  Now, define the penalized criteria by

\begin{equation}\label{Cont_pen}
 \widehat{C}_n(m):= -2 \widehat L_n(\widehat \theta(m)) + \kappa_n |m|, \text{ for all } m \in \mathcal{M} ,
 \end{equation}
 where $(\kappa_n)_{n \in \N}$ is an increasing sequence of the regularization parameter (possibly data-dependent) that will be used to
calibrate the penalty term, 
and
$|m|$ is the number of non-zero components of $\theta^*(m) \in \Theta_m$ that will be called the dimension of the model $m$. 
The selection of the "best" model $\widehat m_n$ is then obtained by minimizing the penalized contrast; that is,
\begin{equation}\label{Estim_m}
  \widehat m_n :=  \underset{  m \in \mathcal{M}}{\text{argmin}} \left(\widehat{C}_n(m)\right).
  \end{equation} 

  \noindent 
 Using the results of Theorems \ref{th1} and \ref{th2}, we establish the asymptotic behavior of the model selection procedure, as shown in the following theorem.
 
 \begin{thm}\label{th4} 
 Let $(Y_1,\ldots,Y_n)$ be a trajectory of a process belonging to $\mathcal {AC}$-$X(M_{\theta^*},f_{\theta^*})$, where  $\theta^* \in \Theta \cap \Theta(r)$ with  $r > 4$. 
Assume that (\textbf{A0})-(\textbf{A2}),  (\textbf{A}$_i(f_\theta,\Theta)$), (\textbf{A}$_i(M_\theta,\Theta)$) (or (\textbf{A}$_i(H_\theta,\Theta)$)) (for $i=0,1,2$) and (\ref{exp_Lip_cov}) (with $r > 4$) hold with $\kappa_n/ n \limiten 0$. Suppose that when $\theta^* \in \overset{\circ}{\Theta}$, 
\begin{multline}\label{cond_th4}
  \sum_{k \geq 1} \frac{1}{\sqrt{k \log \log k}}\sum_{j \geq k} \sum_{i = 0}^{2}
  \big\{ \alpha^{(i)}_{j,Y}(f_\theta,\Theta)+\alpha^{(i)}_{j,X}(f_\theta,\Theta)+\alpha^{(i)}_{j,Y}(M_\theta,\Theta)+\alpha^{(i)}_{j,X}(M_\theta,\Theta)\\
   +\alpha^{(i)}_{j,Y}(H_\theta,\Theta)+\alpha^{(i)}_{j,X}(H_\theta,\Theta)
   \big\} <
   \infty.
\end{multline}
\rm  
\begin{enumerate}
	\item [(i.)] 
	 \it If $\kappa_n/\sqrt{\log \log n} \limiten \infty$, then 
	 \[ \widehat m_n  \limiteproban m^* .\]
	 
	 \rm
 \item [(ii.)]  
 \it When $\theta^* \in \overset{\circ}{\Theta}$, there exists a constant $c$ such that if $ ~\underset{n\rightarrow \infty}{\liminf}  (\kappa_n/\log \log n)>c$, then 
 \[
 \widehat m_n  \limitepsn m^* . 
 \]
 
 \rm
 \item [(iii.)]
 \it If $\theta^* \in \overset{\circ}{\Theta}$ and (\ref{cond_th4}) holds, then 
 \[
 \widehat \theta(\widehat m_n)- \theta^*= O\bigg(\sqrt{\frac{\log \log n}{n}}\bigg).
 \]
\end{enumerate}
 \end{thm}

 \noindent
 Remark that, if  $\sum_{i = 0}^{2}
  \big\{ \alpha^{(i)}_{j,Y}(f_\theta,\Theta)+\alpha^{(i)}_{j,X}(f_\theta,\Theta)+\alpha^{(i)}_{j,Y}(M_\theta,\Theta)+\alpha^{(i)}_{j,X}(M_\theta,\Theta) +\alpha^{(i)}_{j,Y}(H_\theta,\Theta)+\alpha^{(i)}_{j,X}(H_\theta,\Theta) \big\} = O(j^{-\gamma})$ for some $\gamma > 3/2$, then (\ref{cond_th4}) is satisfied. 
 The first and second parts of Theorem \ref{th4} show the consistency of the selection procedure; in particular, the second part provides sufficient conditions for the consistency of the HQC procedure. 
The last part establishes that the estimator of the parameter of the selected model $\widehat \theta(\widehat m_n)$ obeys the law of iterated logarithm.


 \section{Some examples}\label{Sec_Example}
 In this section, we detail some particular processes satisfying the class (\ref{Model}).  
 We show that the regularity conditions required for the main results are satisfied for these processes, with a particular emphasis on the identifiability assumption. 
For each example discussed, we consider that $X_t=(X_{1,t},X_{2,t},\ldots,X_{d_{x}, t}) \in \R^{d_{x}}$ ($d_{x} \in \N$) represents  a vector of covariates; and
   $(\xi_t)_{t \in \Z}$ is a sequence of zero-mean \textit{i.i.d.} random variable satisfying $\E(\xi^r_0) < \infty$ for some $r \geq 2$ and $\E(\xi^2_0) =1$.

 \subsection{Threshold ARX$(\infty)$ models}\label{sub_sect_Tarx}
 Consider the threshold autoregressive model with exogenous covariates (TARX$(\infty)$) defined by, 
\begin{equation}\label{Tarax_model}
 Y_t=  \psi_0(\theta^*) +\sum\limits_{k \geq 1} \Big( \psi_k^+(\theta^*) \max(Y_{t-k},0) + \psi_k^-(\theta^*) \min(Y_{t-k},0) \Big)  + \sum\limits_{k \geq 1} \gamma'_k(\theta^*)X_{t-k} +\xi_{t}, \ \forall t \in \Z,
 \end{equation}
 where $\theta^*$ is the true parameter and $\psi_0(\cdot)$, $\psi_k^+(\cdot)$, $\psi_k^-(\cdot)$, $\gamma_k(\cdot)$ (for  $k \geq 1$) are assumed to be twice  continuously differentiable functions on $\Theta$.
 This model is a generalization of the threshold AR process of Tong (1990). Also, the ARMAX process (see Hannan and Deistler (2012)) is a  specific example of the model (\ref{Tarax_model}).  
%
 Set $f^t_\theta = \psi_0(\theta) +\sum\limits_{k \geq 1} \Big( \psi_k^+(\theta) \max(Y_{t-k},0) + \psi_k^-(\theta) \min(Y_{t-k},0) \Big)  + \sum\limits_{k \geq 1} \gamma'_k(\theta)X_{t-k}$ for all $\theta\in\Theta$.  
If $\sum\limits_{k \geq 1} \| \gamma'_k (\theta^*) \| < \infty $ and $ \sum\limits_{k \geq 1}\max \left\{\alpha_k(g), |\psi_k^+(\theta^*)|, |\psi_k^-(\theta^*)| \right\}<1$, then, assumption \textbf{A}$_0 (f_\theta,\{ \theta^*\})$ holds and there exists a stationary and ergodic solution with $r$-order moment. The assumption (\textbf{A1}) holds with $\underline{h}=1$. 
 Denote for $t, i \in \Z$,
  \begin{equation}\label{F_t_i}
  \mathcal F_{t,i} = \sigma( \xi_{t-j}, j > i, \ X_{t-k}, k>0 ),
\end{equation}   
the $\sigma-$field generated by $ \left\{ \xi_{t-j}, j > i, \ X_{t-k}, k>0\right\}$.
Let us set following  additional conditions:
 \begin{enumerate}
  \item [(\textbf{B0}):] $\E \left[ \xi_t X_{t'}\right]=0$ for all $(t,t') \in \Z^2$;
  \item [(\textbf{B1}):] for $c^+, c^- \in \R$ such that $c^+ \neq 0 $ or $ c^- \neq 0$, $ c^+ \max(Y_{t-i},0) + c^- \min(Y_{t-i},0)$ given $\mathcal F_{t,i}$ is non-degenerate;  
  \item [(\textbf{B2}):] if $(c_k)_{k \in \N}$ is a sequence  of vector of $\R^{d_x}$ such as  $\exists  c_{k_0} \neq 0$ (with $k_0 \in \N$), then 
$\sum\limits_{k \geq 1} c'_kX_{t-k}$ is non-degenerate;
  \item [(\textbf{B3}):] the function $\theta \mapsto  \psi_{k_0}^+(\theta)$ or $\theta \mapsto  \psi_{k_0}^-(\theta)$, for some $k_0 \geq 1$ is injective; 
 \item [(\textbf{B4}):] The function $\theta \mapsto \psi_0(\theta)$ (or $\theta \mapsto \gamma_{k_0}(\theta)$ for some $k_0 > 0$) is  is injective  and  holds (\textbf{B2}).
  \end{enumerate}
  
  \medskip
  
  %
  To ensure the identifiability, both the assumptions (\textbf{B3}) and (\textbf{B4}) are not necessary;  that is, the model is identifiable
if (\textbf{B1}) and ((\textbf{B3}) or (\textbf{B4})) hold.
%
%
 Indeed, let $\theta\in\Theta$ such that $f^t_\theta = f^t_{\theta^*} $. Then,
\begin{multline} \label{Tarx_eq_ft_theta}
\sum\limits_{k \geq 1} \Big( \big(\psi_k^+(\theta^*) - \psi_k^+(\theta) \big) \max(Y_{t-k},0) + \big(\psi_k^-(\theta^*) - \psi_k^-(\theta) \big) \min(Y_{t-k},0) \Big)  \\ 
= \psi_0(\theta) - \psi_0(\theta^*) +  \sum\limits_{k \geq 1} \big(\gamma'_k(\theta) - \gamma'_k(\theta^*) \big)X_{t-k}. 
\end{multline} 
 By contradiction, assume that  $\psi_k^+(\theta^*) - \psi_k^+(\theta) \neq 0$ or $\psi_k^-(\theta^*) - \psi_k^-(\theta) \neq 0$ for some $k >0$ and let $m > 0$ be the smallest integer satisfying  $\psi_m^+(\theta^*) - \psi_m^+(\theta) \neq 0$ or $\psi_m^-(\theta^*) - \psi_m^-(\theta) \neq 0$.  
 It holds from (\ref{Tarx_eq_ft_theta}) that,
 \begin{align} \label{Tarx_varfi_ft_theta} 
 \nonumber &\big(\psi_m^+(\theta^*) - \psi_m^+(\theta) \big) \max(Y_{t-m},0) + \big(\psi_m^-(\theta^*) - \psi_m^-(\theta) \big) \min(Y_{t-m},0) \\
 \nonumber &=   \psi_0(\theta) - \psi_0(\theta^*) 
   - \sum\limits_{k > m} \Big( \big(\psi_k^+(\theta^*) - \psi_k^+(\theta) \big) \max(Y_{t-k},0) + \big(\psi_k^-(\theta^*) - \psi_k^-(\theta) \big) \min(Y_{t-k},0) \Big) \\ 
   & \hspace{10cm} +  \sum\limits_{k \geq 1} \big(\gamma'_k(\theta) - \gamma'_k(\theta^*) \big)X_{t-k}. 
\end{align} 
Since the right-hand side of (\ref{Tarx_varfi_ft_theta}) is $\mathcal F_{t,m}$-measurable (thanks to the non anticipative property of the process $(Y_t)$), $\big(\psi_m^+(\theta^*) - \psi_m^+(\theta) \big) \max(Y_{t-m},0) + \big(\psi_m^-(\theta^*) - \psi_m^-(\theta) \big) \min(Y_{t-m},0)$ given $\mathcal F_{t,m}$ is degenerate, which contradicts the assumption $(\textbf{B1})$. Thus,  $\psi_k^+(\theta^*) - \psi_k^+(\theta) = 0$ and $\psi_k^-(\theta^*) - \psi_k^-(\theta) =0$ for all $k \geq 1$. Therefore,
\begin{itemize}
	\item   if (\textbf{B3}) holds, then $\theta =\theta^*$;
	\item else, if (\textbf{B4}) holds, then we have from (\ref{Tarx_varfi_ft_theta}),
	\[
 \sum\limits_{k=1}^{\infty} (\gamma'_k(\theta^*)-\gamma'_k(\theta))X_{t-k} = 
\psi_0(\theta)-\psi_0(\theta^*), 
	\]
which implies that  $\gamma_k(\theta^*)=\gamma_k(\theta)\ \forall k\geq 1$ (by the assumption (\textbf{B2})); and consequently,  $\psi_0(\theta)=\psi_0(\theta^*)$.  Hence, $\theta=\theta^*$.  
\end{itemize}
In the case of the ARX($\infty$) models (obtained when $\psi_k^+(\theta^*) = \psi_k^-(\theta^*)$, for all $k\geq1$), the condition (\textbf{B1}) is not necessary; the identifiability holds with (\textbf{B0}) and ((\textbf{B3}) or (\textbf{B4})).   

\subsection{Asymmetric Power ARCH-X$(\delta,\infty)$ models} \label{sub_sect_APARCHX_infty}
  Consider the Asymmetric Power ARCH with with exogenous covariates (APARCH-X$(\delta,\infty)$) defined by,
 \begin{align}\label{APARCHX_infty_model}
 & \nonumber Y_t=  \xi_{t} \sigma_t((\theta^*)
  \text{ with } \\
 & \sigma^\delta_t(\theta^*)=\phi_0(\theta^*) +\sum\limits_{k=1}^{\infty} \Big( \phi_k^+(\theta^*) \big( \max(Y_{t-k},0) \big)^\delta  + \phi_k^-(\theta^*)\big( \max(-Y_{t-k},0) \big)^\delta  \Big) +  \sum\limits_{k=1}^{\infty} \gamma'_k(\theta^*)X_{t-k}, \ \forall t \in \Z,
 \end{align}
 where  $\theta^*$ is the true parameter, $\phi_0(\cdot)$, $\phi_k^+(\cdot)$, $\phi_k^-(\cdot)$, $\gamma_k(\cdot)$, $k\geq 1$ are non-negative (componentwise  for $\gamma_k(\theta^*)$) functions assumed to be twice  continuously differentiable on $\Theta$,  with $(\phi_0(\theta)) > 0$ for all $\theta\in\Theta$, $X_t$ is a vector of non-negative (componentwise) covariates, and $\delta>0$. This process is an example of the class (\ref{Model}) with $f^t_\theta=0$ and $M^t_\theta=\sigma_t(\theta)$. 
 Numerous classical ARCH-type parametrizations, for instance, GARCH-X (obtained with $\delta = 2$ and $\phi_k^+(\theta^*) = \phi_k^-(\theta^*)$), TARCH-X (obtained with $\delta = 1$) are specific example of (\ref{APARCHX_infty_model}).
 There are several works in the literature based on the GARCH-X model, see for instance Han and  Kristensen (2014), Nana {\it et al.} (2013), Han (2015). 
 Model (\ref{APARCHX_infty_model}) is a generalization of the class of APARCH-X($\delta, p,q$) studied by Francq and Thieu (2019).
 Assume that $\theta^* \in \Theta(r)$, where $\Theta(r)$ is given by
%
 \[
\Theta(r)= \big\{ \theta \in \R^{d} ~ \big / ~ \sum\limits_{k \geq 1}\max \left\{\alpha_k(g), \| \xi_0 \|_r |\phi_k^+(\theta)|^{1/\delta}, \| \xi_0 \|_r  |\phi_k^-(\theta)|^{1/\delta} \right\}<1
\text{ and } \sum\limits_{k \geq 1}\|\gamma_k(\theta)\|^{1/\delta} <\gamma_U \text{ for some } \gamma_U >0 \big\}.
\]
Therefore, \textbf{A}$_0 (H_\theta,\{ \theta^* \})$ holds; and a stationary and ergodic solution with $r$-order moment exists. If  $\inf_{\theta \in \Theta}\phi_0(\theta)>0$, then the assumption (\textbf{A1}) is satisfied. 
The following assumptions are needed to ensure the identifiability.
  \begin{enumerate}
 
 \item [(\textbf{B5}):]  for all $i\geq1$ and $t \in \Z$, the support of the distribution of $\xi_{t-i}$ given $\mathcal F_{t,i}$ is not included in $[0 ,\infty)$ or in $(-\infty, 0]$ and contains at least three points.
  
  \item [(\textbf{B6}):] the function $\theta \mapsto \phi_{k_0}^+(\theta)$ or $\theta \mapsto \phi_{k_0}^-(\theta)$ for some $ k_0 > 0$, is injective. 
  
  \item [(\textbf{B7}):]  The function $\theta \mapsto \phi_{0}(\theta)$  (or $\theta \mapsto \gamma_{k_0}(\theta^*)$ for some $ k_0 > 0$) is injective, and the condition (\textbf{B2}) holds for this model. 
  \end{enumerate}
By going along similar lines as in the Subsection \ref{sub_sect_Tarx} and by using the Lemma 4 in \cite{Francq2019}, one get that, if (\textbf{B5}) and ((\textbf{B6}) or (\textbf{B7})) hold, then the model (\ref{APARCHX_infty_model}) is identifiable. 

\medskip
%
%
Let us stress that, the {\it i.i.d.} assumption for $(\xi_t)_{t \in \Z}$ is a bit strong for the model (\ref{APARCHX_infty_model}). 
This assumption, which is needed to the large class $\mathcal {AC}$-$X(M_{\theta^*},f_{\theta^*})$ can be relaxed to $(\xi_t, \mathcal F_{t,0})_{t \in \Z}$ is a martingale difference sequence (see, for instance, \cite{Francq2019} in the case of APARCH-X($\delta, p,q$) model) when checking the identifiability.
Also, as pointed out by Francq and Thieu (2019), in the absence of covariates and when $(\xi_t)$ is \textit{i.i.d} (i.e., the case of the APARCH$(\infty)$ model), the assumption (\textbf{B5}) can be automatically reduced to: $P(\xi_0 >0) \in (0,1)$ and the support of the distribution of $\xi_0$ contains at least three points. 
Moreover, Assumption (\textbf{B5})  prevents taking redundant covariate; for instance, it excludes the situation where $X_{t-1}=\max(Y_{t-i},0)^\delta$ or $\max(-Y_{t-i},0)^\delta$ for some $i\geq 1$. 


\subsection{ARX($\infty$)-ARCH($\infty$) models}
 Consider the ARX($\infty$)-ARCH($\infty$) model given by, 
  \begin{equation}\label{ARX_ARCH_infinite_model} 
\left\{
\begin{array}{l}
Y_t= \psi_0(\theta^*) +\sum\limits_{k \geq 1} \psi_k(\theta^*) Y_{t-k} + \sum\limits_{k \geq 1} \gamma'_k(\theta^*)X_{t-k} +\varepsilon_t\\
\rule[0cm]{0cm}{.5cm}
\varepsilon_t =\xi_t \sigma_t(\theta^*) ~\text{ with }~ \sigma^2_t(\theta^*)=\phi_0(\theta^*) +\sum\limits_{k \geq1} \phi_k(\theta^*) \varepsilon^2_{t-k},
\end{array}
\right.
   \end{equation} 
 where $\theta^*$ is the true parameter, and $\psi_0(\cdot)$, $\psi_k(\cdot)$, $\gamma_k(\cdot)$, $\phi_0(\cdot)$, $\phi_k(\cdot)$, $k \geq 1$ are  assumed to be twice  continuously differentiable on $\Theta$, and satisfying $\phi_0(\theta) > 0$ for all $\theta \in \Theta$.
 Model (\ref{ARX_ARCH_infinite_model}) is an extension of the ARMA-GARCH, ARMAX-GARCH processes.
 This model belongs to the class $\mathcal {AC}$-$X(M_{\theta^*},f_{\theta^*})$ with 
 \begin{align*}
 f^t_{\theta}&=\psi_0(\theta) +\sum\limits_{k \geq 1} \psi_k(\theta) Y_{t-k} + \sum\limits_{k \geq 1} \gamma'_k(\theta)X_{t-k}  ~ ~ \text{ and }\\
 M^t_{\theta}&=\sqrt{ \phi_0(\theta) +\sum\limits_{k \geq 1} \phi_k(\theta) \Big( Y_{t-k} - \psi_0(\theta) -\sum\limits_{j \geq 1} \psi_j(\theta) Y_{t-k-j} - \sum\limits_{j \geq 1} \gamma'_j(\theta)X_{t-k-j} \Big)^2 }, 
\end{align*}
 for all $\theta \in \Theta$. Hence, the assumption \textbf{A}$_0 (f_\theta,\{ \theta\})$ holds with
 $\alpha^{(0)}_{k,Y}(f_{\theta},\{\theta\}) = | \psi_k(\theta) | $ and $\alpha^{(0)}_{k,X}(f_{\theta},\{\theta\}) = \| \gamma_k(\theta) \| $.
 From an expansion of $M^t_{\theta}$, one can easily get that, \textbf{A}$_0 (M_\theta,\{ \theta\})$ holds with
  \begin{equation*}
\left\{
\begin{array}{l}
\alpha^{(0)}_{1,Y}(M_{\theta},\{\theta\})= \sqrt{ \phi_k(\theta) }  ~\text{ and }~ \alpha^{(0)}_{k,Y}(M_{\theta},\{\theta\})= \sqrt{ \phi_k(\theta) } + \sum_{i=1}^{k-1} \sqrt{ \phi_i(\theta) } | \psi_{k-i}(\theta) | ~~ \text{ for } k \geq 2;\\
\rule[0cm]{0cm}{.6cm}

\alpha^{(0)}_{1,X}(M_{\theta},\{\theta\})=0~ \text{ and }~\alpha^{(0)}_{k,X}(M_\theta,\{\theta\})=  \sum_{i=1}^{k-1} \sqrt{ \phi_i(\theta) } \| \gamma_{k-i}(\theta) \| ~~ \text{ for } k \geq 2.
\end{array}
\right.
   \end{equation*} 
Therefore,  the stationarity set $\Theta(r)$ is defined by
 \begin{multline*}
\Theta(r)= \Big\{ \theta \in \R^{d} ~ \big / ~   \sum\limits_{k \geq 1}\max \Big\{   \alpha_k(g), | \psi_k(\theta) | + \| \xi_0 \|_r \Big( \sqrt{ \phi_k(\theta) } + \sum_{i=1}^{k-1} \sqrt{ \phi_i(\theta) } | \psi_{k-i}(\theta) | \Big)  \Big\}<1
\text{ and } \\
\sum\limits_{k \geq 1}\|\gamma_k(\theta)\|  <\gamma_U  \text{ for some } \gamma_U >0 \Big\}.
\end{multline*} 
 Assumption (\textbf{A1}) is satisfied if $\inf_{\theta \in \Theta}\phi_0(\theta)>0$. The identifiability conditions can be obtained as in Subsection \ref{sub_sect_Tarx} and \ref{sub_sect_APARCHX_infty}.


\section{Simulation study }
In this section, we consider a double autoregressive model with exogenous covariates, defined by   
  \begin{equation}\label{FDAR_model} 
Y_t= \phi_0 + \sum_{i=1}^{p_1} \phi_i  Y_{t-i} +\sum_{i=1}^{q_1} \psi_i' X_{t-i}+\xi_t \sqrt{\alpha_0 + \sum_{i=1}^{p_2} \alpha_i Y^2_{t-i} +\sum^{q_2}_{i=1} \beta_i' ( X_{t-i} \odot X_{t-i}) },
   \end{equation}  
 where $(X_t)_{t \in \Z}$ is an exogenous multivariate covariate process with values in $\R^{d_x}$ ($d_x \in \N$),  
 $p_1, p_2, q_1, q_2 \in \N$,  $\phi_0, \phi_1,\cdots, \phi_{p_1} \in \R$, $\psi_1,\cdots,\psi_{q_1} \in \R^{dx}$,
  $\alpha_0>0$, $\alpha_1, \cdots, \alpha_{p_2} \geq0$, $\beta_1,\cdots,\beta_{q_2} \in [0,\infty)^{d_x}$, $\odot$ denotes the Hadamard product (componentwise multiplication),  $\xi_t$ is a white noise with $\E \xi_0^2 = 1$.
 This model is a generalization of the factor double autoregressive (FDAR) process introduced by Guo \textit{et al.} (2014) to extend the double AR$(p)$ model proposed by Ling (2007). The ARX($p$) and ARCH-X($p$) are particular cases of the model (\ref{FDAR_model}).
 We assume that  $(X_t)_{t \in \Z}$ is a VAR(1) (vector autoregressive) process:  
 \begin{equation}\label{VAR1}
       X_{t}=\varphi_0 +\varphi_1 X_{t-1} + \eta_{t}~ \text{ for all }  t \in \Z,
 \end{equation}
  where $\varphi_0 \in \R^{d_x}$, $\varphi_1$ is a real coefficients $(d_x \times d_x)$-matrix and $\eta_{t}$  is a white noise with $\E( \eta_{0} \eta_{0}') = \Sigma_{\eta}$. 
 If $\| \varphi_1\| < 1$, then the stability condition (\ref{exp_Lip_cov}) holds with $\alpha_1(g) = \| \varphi_1\|$ and $\alpha_k(g)=0$ for $k \geq 2$. 
 The stationarity set $\Theta(r)$ is
 \begin{multline*}
  \Theta(r) = \big\{
   \theta=(\phi_0,\phi_1,\cdots,\phi_{p_1}, \alpha_0,\alpha_1,\cdots, \alpha_{p_2}, \psi_1',\cdots,\psi_{q_1}',\beta_1',\cdots,\beta'_{q_2} ) \in \Theta, \text{ where } \\
   \Theta \subset \R^{p_1+1} \times (0,\infty) \times [0,\infty)^{p_2} \times \R^{d_x q_1} \times [0,\infty)^{d_xq_2} 
   \text{ and }
      \sum_{i=1}^{\max(1,p_1, p_2)} \max \left\{\alpha_i(g),|\phi_i|+ \|\xi_0\|_r \alpha_i \right\}<1 \big\},
  \end{multline*}
  with $\phi_i=0$ if $i>p_1$ and $\alpha_i=0$ if $i>p_2$.
 Based on the examples discussed in Section \ref{Sec_Example}, if the conditions (\textbf{B0}) and (\textbf{B3}) hold for (\ref{FDAR_model}), then to satisfy the identifiability condition, it suffices to impose  the following assumption on the covariate:

  \medskip 
  
\noindent (\textbf{B7}): if $(c_{k,1})_{1\leq k \leq q_1}$ and $(c_{k,2})_{1\leq k \leq q_2}$ are sequences  of vector of $\R^{d_x}$ such as  $\exists  c_{k_1,1} \neq 0$ and $\exists  c_{k_2,2} \neq 0$ (with $1\leq k_1 \leq q_1$, $1\leq k_2 \leq q_2$), then $\sum^{q_1}_{i=1} c_{k,1}' X_{t-k}$ and  $\sum^{q_2}_{i=1} c_{k,2}' X_{t-k} \odot X_{t-k}$ are not degenerated.  
 
 \medskip
 
 Set $ \psi=( \psi_1,\ldots, \psi_q)$ and $ \beta=( \beta_1,\ldots, \beta_q)$; the true parameter is $\theta^*=(\phi_0,\phi_1,\alpha_0,\alpha_1,\psi,\beta)$. In the sequel, we focus on the following two cases.
 
 \medskip
 \noindent \textbf{Case 1}.
  We consider an example of model (\ref{FDAR_model}) with univariate covariates where $p_1=p_2=1$ and $q_1=q_2=q$; the AR parameter is set to $(\varphi_0,\varphi_1)=(0.5,0.5)$ and $\Sigma_\eta=1$.
%

\medskip
 \noindent \textbf{Case 2}. In this second example, model (\ref{FDAR_model}) is considered with univariate/multivariate covariates where $p_1=p_2=1$, $q_1=1$ and $q_2=0$; 
 thus, the true parameter is $\theta^*=(\phi_0,\phi_1,\alpha_0,\alpha_1,\psi_1')$.
   This second example is related to the real data application (see Section \ref{sect_Real_data}). 
%

\subsection{Estimation and significance test}
  Some results from Monte Carlo simulations are displayed to assess the asymptotic properties of the QMLE, as well as 
 to investigate the empirical size and power of the proposed procedure on testing the significance of the covariate $X_t$.
We will consider samples where the innovation $(\xi_t)_{t \in \Z}$ is generated from Gaussian and Student distributions (with 5 with degrees of freedom).  
 The model (\ref{FDAR_model}) is considered in the \textbf{Case 1} with $p_1=p_2=q_1=q_2=1$  (scenario {\bf S}$_0$, {\bf S}$_1$, {\bf S'}$_0$ and {\bf S'}$_1$ below).
  In the \textbf{Case 2} ($p_1=p_2=q_1=1$, $q_2=0$), we consider scenarios with univariate covariate (scenario {\bf S''}$_0$ and {\bf S''}$_1$ below), the AR parameter are  $(\varphi_0,\varphi_1,\Sigma_\eta)=(23.61, 0.7,21.56)$.
%
%
\begin{itemize}
    \item {\bf scenario S}$_0$: $\theta^*=(0.15,-0.2,0.4,0.3,0,0)$;
    \item {\bf scenario S}$_1$: $\theta^*=(0.15,-0.2,0.4,0.3,0.08,0)$;
     \item {\bf scenario S'}$_0$: $\theta^*=(1,0.4,0.5,0.2,0,0)$;
    \item {\bf scenario S'}$_1$: $\theta^*=(1,0.4,0.5,0.2,0.07,0.07)$.
    \item {\bf scenario S''}$_0$: $\theta^*=(37.95,0.33, 32.11, 0.02,0)$;
    \item {\bf scenario S''}$_1$: $\theta^*=(37.95, 0.33, 32.11, 0.02, -0.21)$. 
\end{itemize}
The scenarios \textbf{S}$_0$, \textbf{S'}$_0$ and \textbf{S''}$_0$ correspond to cases where the covariate is absent; \textbf{S''}$_0$ and \textbf{S''}$_1$ are related to the real data application.
We consider the following significance tests:
\begin{equation*}
 H_0: \, \theta^*=(0.15,-0.2,0.4,0.3,0,0) ~ ~ (S_0) ~~~~~~ \text { against } ~~~~~~  H_1: \,  \theta^* \neq (0.15,-0.2,0.4,0.3,0,0);
\end{equation*}
 \medskip
\begin{equation*}
 H_0: \, \theta^*=(1,0.4,0.5,0.2,0,0) ~ ~ (S'_0) ~~~~~~ \text { against } ~~~~~~  H_1: \,  \theta^* \neq (1,0.4,0.5,0.2,0,0).
\end{equation*}
 \medskip
\begin{equation*}
 H_0: \, \theta^*=(37.95,0.33, 32.11, 0.02,0) ~ ~ (S''_0) ~~~~~~ \text { against } ~~~~~~  H_1: \,  \theta^* \neq (37.95,0.33, 32.11, 0.02,0).
\end{equation*}
 
\medskip

\begin{figure}[h!]
\begin{center}
\includegraphics[height=8cm, width=15cm]{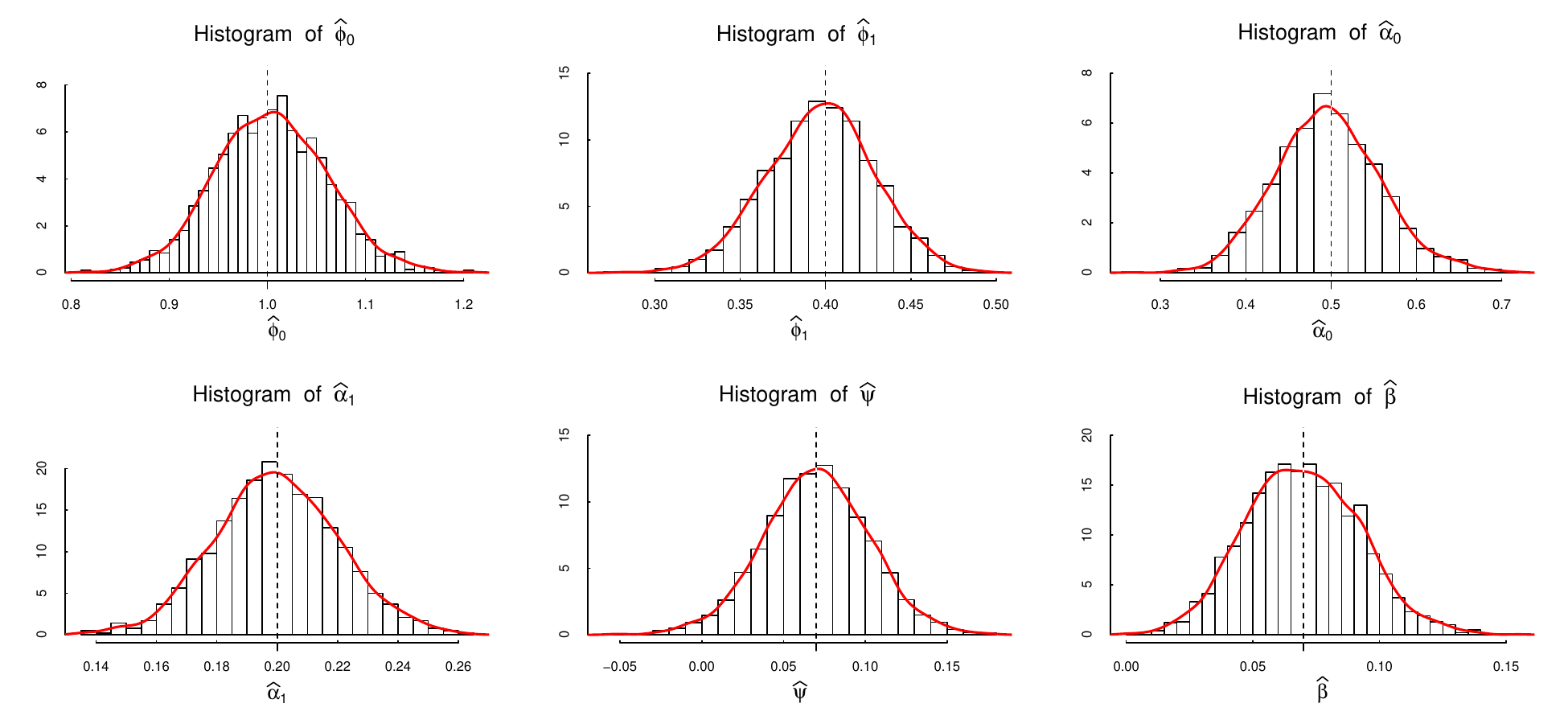}
\end{center}
\vspace{-.6cm}
\caption{\it Histograms of the components of $\widehat \theta_n$ in the scenario \textbf{S'}$_1$ with sample size $n=1000$.
 The overlaying curves are the density estimates and the dotted vertical lines represent the true values of the parameters.}
\label{Graphe_normal}
\end{figure}

\noindent In each of the scenarios \textbf{S}$_0$, \textbf{S}$_1$, \textbf{S'}$_0$, \textbf{S'}$_1$, \textbf{S''}$_0$ and \textbf{S''}$_1$, we simulate $200$ replications with the sample size $n = 500,\,1000$ and test the nullity of the vector $(\psi,\beta)$ after estimating the parameters of interest.
Table \ref{Table_Res} contains the empirical mean and root mean square error (RMSE) of each component of the estimator.
The last two columns of Table \ref{Table_Res} indicate the empirical levels and powers of the above tests at the nominal level $\alpha=0.05$, where the empirical powers are computed
under the alternative $H_1$ respectively in the scenario \textbf{S}$_1$, \textbf{S'}$_1$ and  \textbf{S''}$_1$.
%
 For the scenario {\bf S'}$_1$, the histograms and estimated densities of the estimates are plotted in Figure \ref{Graphe_normal}. 

\medskip

From these findings, one can see that, in all scenarios, the performance of the QMLE is satisfactory in terms of the mean and that, the RMSE of the estimators decreases when $n$ increases. This is consistent with the results of Theorem \ref{th1}. 
Also remark that, the fact of computing the QMLE with $(\widehat \psi, \widehat \beta)$ for trajectories generated without covariates (see the scenarios \textbf{S}$_0$, \textbf{S'}$_0$ and  \textbf{S''}$_0$) does not affect the performance of the QMLE, which again confirms its good theoretical properties. 
As seen in Figure \ref{Graphe_normal}, for each component of $\widehat \theta_n$, the estimated density is very close to that of the normal distribution; which is in accordance with the asymptotic results obtained from Theorem \ref{th2} when $\theta^* \in \overset{\circ}{\Theta}$.
The results of the test (see Table \ref{Table_Res}) show that, the statistic $W_n$ is slightly oversized for  $n=500$ in cases where the innovation
 is generated from Student distributions, but the empirical levels are reasonable when $n = 1000$ in the sense that, they are very close  to the nominal one.
Further, the empirical powers displayed increases with the sample size and are quite accurate. 

%

\subsection{Model selection}
Now, we are going to carry out other simulation experiments aimed at evaluating the effectiveness of the proposed model selection procedure in the model (\ref{FDAR_model}) for choosing the order $q_1=q_2=q$ in the \textbf{Case 1}.
To this end, $q=2$ is set as the "true" model $m^*$ and that the following scenarios are considered: 
\begin{itemize}
    \item {\bf scenario S}$^*_1$: $\theta^*=(0.6,0.45,0.5,0.15, 1,0.7, 0.6,0.35)$;
    
    \item {\bf scenario S}$^*_2$: $\theta^*=(0.15,0.4,0.5,0.2, 0.1,0.1, 0.03,0.3)$.
\end{itemize}
The competing models used are all process satisfying (\ref{FDAR_model}) with $q \in \{0,1,\ldots,9 \}$, which leads us to a collection of 10 models.

\medskip 
  
  \noindent
 In the \textbf{Case 2}, consider the multivariate covariate $X_t=(X_{t,1}, X_{t,2},\cdots,X_{t,5})$, a VAR(1) (see (\ref{VAR1})) with parameter $\varphi_0 = (23.61, 4.95, 12.12, 716.70, 13.01)'$, $\varphi_1= diag(0.7, 0.57, 0.51, 0.26, 0.56) $ and
\begin{equation*}
\Sigma_\eta = 
\begin{pmatrix}
42.62 & -2.59 &  11.14  & 281.05 & -14.65  \\
-2.59 &  7.50 &  -1.29 &  200.64 &  -8.24 \\
 11.14 & -1.29 &  41.13 & 1019.57  & -1.10  \\
281.05 &  200.64 & 1019.57 & 64948.60 & -745.27 \\ 
-14.64 &  -8.24  & -1.10 & -745.27 & 59.00
\end{pmatrix}
\end{equation*}
Consider the scenario {\bf  S''}$_1$ with the covariate $X_{t,1}$ as the true model.
The \textbf{Case 2} with all the combination of the covariates is performed on the data; that is, there are 32 competing models. This example is related and close to the real data  application. 
 
\begin{figure}[h!]
\begin{center}
\includegraphics[height=12cm, width=16cm]{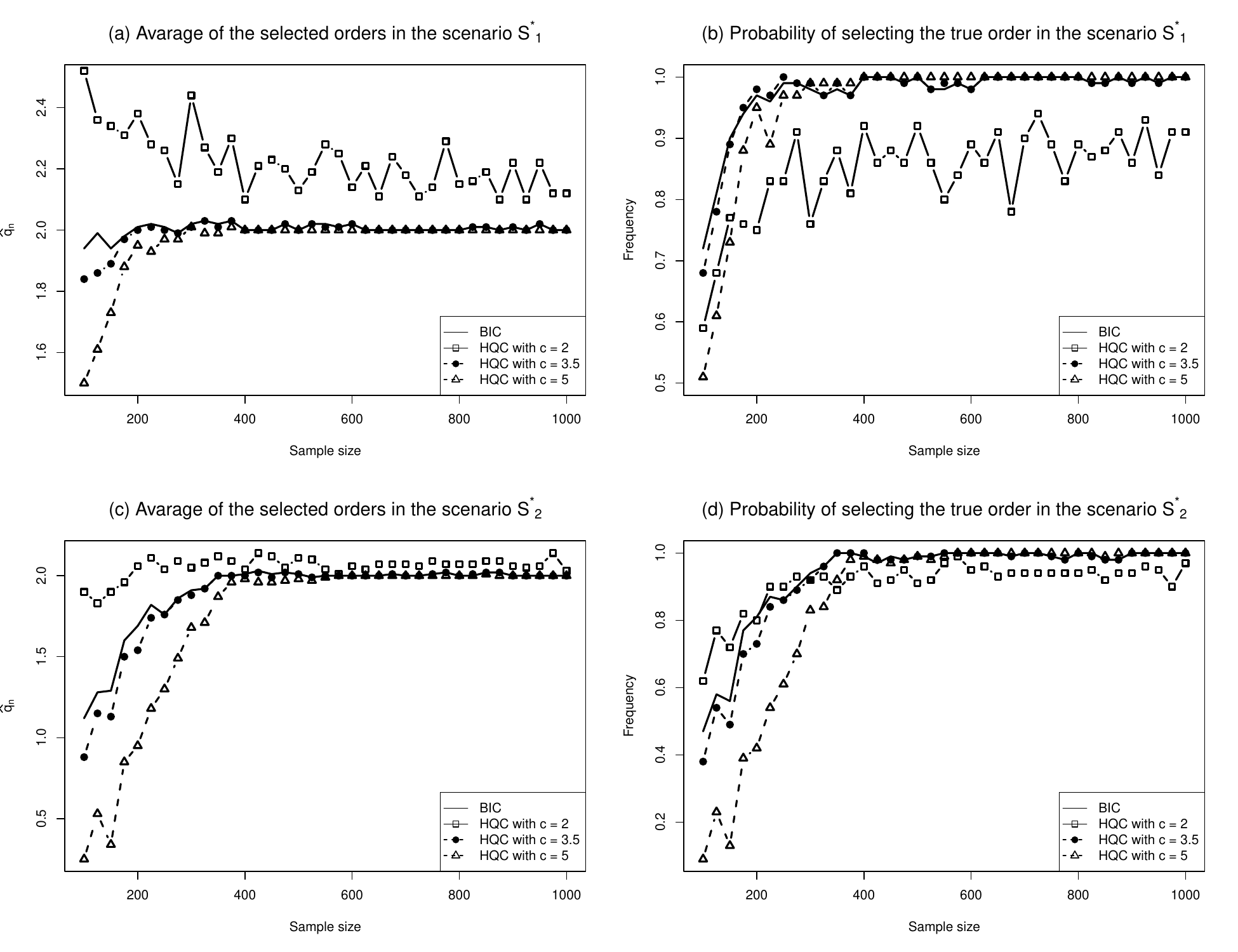}
\end{center}
\vspace{-.6cm}
\caption{\it The averages of the orders selected and the frequencies of selecting the true order based on 100 independent replications depending on sample's length in the scenarios ${\bf S}^*_1$ and ${\bf S}^*_2$.}
\label{Graphe_Selec_Case1}
\end{figure}

\begin{figure}[h!]
\begin{center}
\includegraphics[height=8cm, width=12cm]{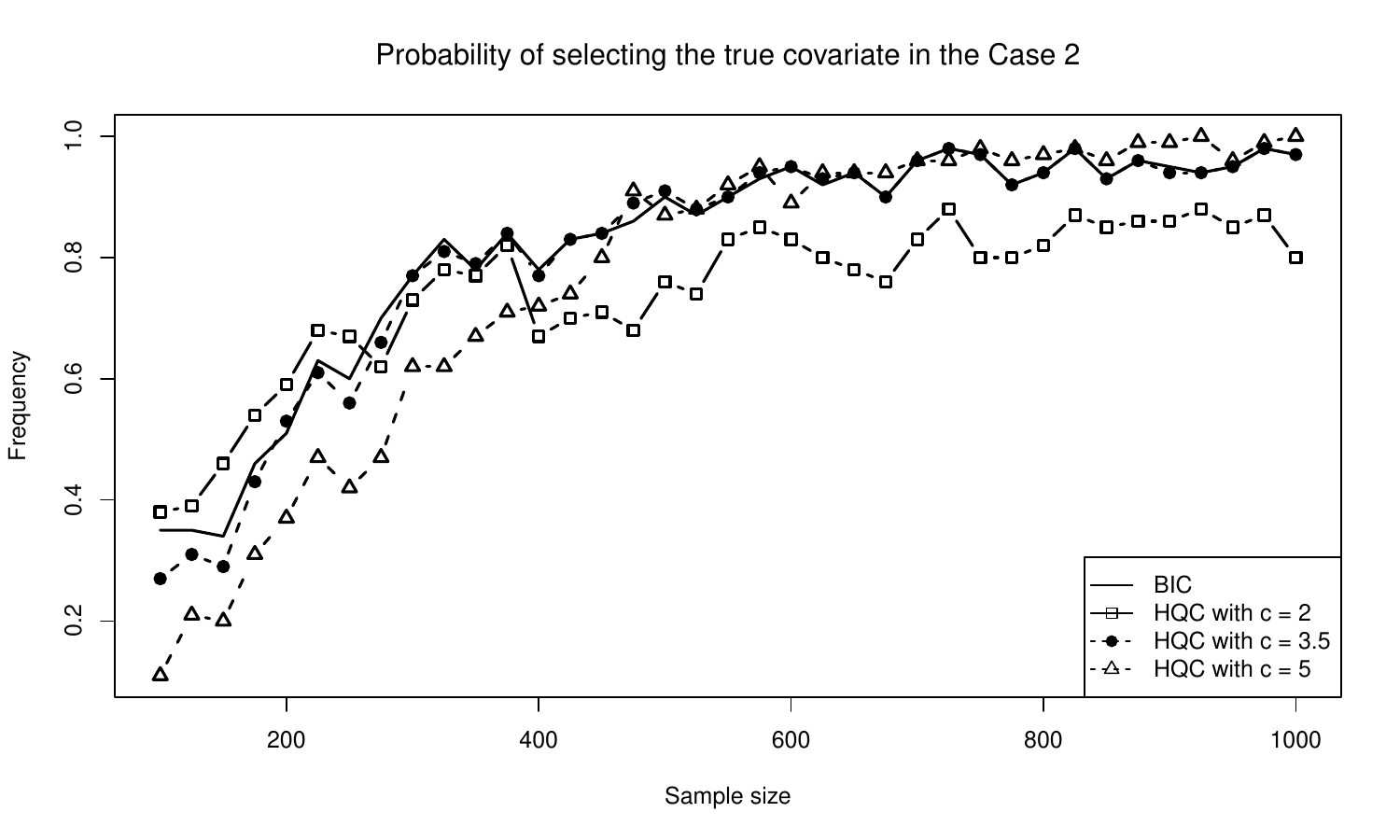}
\end{center}
\vspace{-.6cm}
\caption{\it The frequencies of selecting the true covariate based on 100 independent replications depending on sample's length in the scenario
{\bf S''}$_1$.}
\label{Graphe_Selec_Case2}
\end{figure}

\medskip
\noindent
For $n=100, 125,150,\ldots,1000$, we simulate $100$ independent replications in each of the three scenarios {\bf S}$^*_1$, {\bf S}$^*_2$ and {\bf  S''}$_1$.
We compare the performances of the procedure with $\kappa_n = \log n$ (see  (\ref{Cont_pen})) linked to the Bayesian Information Criteria (BIC) and the procedure with $ \kappa_n = c \log \log n$  ($c \in \{ 2,\, 3.5,\,5\}$) linked to the Hannan-Quinn information Criterion (HQC). 
For the scenarios {\bf S}$^*_1$ and {\bf S}$^*_2$, Figures \ref{Graphe_Selec_Case1} displays the points $(n,\widehat q_n)$, where $\widehat q_n$ denotes the average of the orders selected with trajectories of length $n$, as well as the curve of the proportions of number of replications (frequencies) where the associated criterion selects the true order.
For the scenario of the \textbf{Case 2} (i.e, {\bf  S''}$_1$), the probabilities of choosing the true covariate are displayed in Figure \ref{Graphe_Selec_Case2}. 

\medskip
\noindent
From these figures, 
the first remark is that, for all the penalties, the performances of the procedure increase with $n$ in each scenario. 
Further, the  probability of selecting the true order is very close to $1$ when $n=1000$. This shows that, these procedures are in accordance with the results of Theorem \ref{th4}. 
One can notice that, 
in the scenario ${\bf S}^*_1$, the $\log n$ penalty is more interesting for selecting the true order than the others penalties for a small sample size (see Figure  \ref{Graphe_Selec_Case1} ((a) and (b)) for $n \leq 250$), while in the scenarios ${\bf S}^*_2$ and \textbf{S''}$_1$, the HQC with $c=2$ slightly outperforms the BIC penalization when $n \leq 350$ (see Figure \ref{Graphe_Selec_Case1} ((c) and (d)) and  Figure \ref{Graphe_Selec_Case2}).
However, 
the larger the sample size, the $c\log \log n$ penalty (except for the case where $c=2$) provides the same accuracies in comparison with the $\log n$ penalty, and displays satisfactory results.
The results also show that, as $c$ increases, the performances of the $c\log \log n$ penalty increase, which reveals that the common use of the classical HQC penalization (i.e, the $c\log \log n$ penalty with $c=2$) is not always the optimal choice to select the best model with this information criterion.
%

\begin{table}
\scriptsize
\centering
\caption{\it 
Sample mean and RMSE of the QMLE for the model (\ref{FDAR_model}) following the scenarios  \textbf{S}$_0$, \textbf{S}$_1$, \textbf{S'}$_0$,  \textbf{S'}$_1$, \textbf{S''}$_0$ and \textbf{S''}$_1$. The last two columns show the empirical levels and powers at the nominal level $0.05$ for the test of the relevance of the exogenous covariates.
}
\label{Table_Res}
\vspace{.2cm}
\begin{tabular}{ccccccccccccc}

\Xhline{.7pt}
\rule[0cm]{0cm}{.4cm}
& && \multicolumn{6}{l} {QMLE} &&&  \multicolumn{2} {l}{Statistic $W_n$}   \\
\cline{4-10} \cline{12-13} 
\rule[0cm]{0cm}{.35cm}
Scenario&Noise&$n$&&$\widehat \phi_0$&$\widehat \phi_1$&$\widehat \alpha_0$ &$\widehat \alpha_1$&$\widehat \psi$&$\widehat \beta$
&& Levels& Powers\\
\Xhline{.7pt}
\rule[0cm]{0cm}{.3cm}
  \rule[0cm]{0cm}{.25cm}
   {\bf S}$_0$        &Gaussian& $500$& Mean         &$0.1488$&$-0.2001$&$0.3887$&$0.2987$&$0.0019$&$0.0033$
 &&\multirow{2}*{$0.035$}&\\
\rule[0cm]{0cm}{.25cm}
            &&        & Rmse         &$0.0435$&$0.0607$&$0.0421$&$0.0924$&$0.0276$&$0.0063$ &&&\\

           \rule[0cm]{0cm}{.4cm}
            && $1000$& Mean         &$0.1502$&$-0.2023$&$0.3911$&$0.2979$&$-0.0006$&$0.0032$ 
            &&\multirow{2}*{$0.045$}&\\
\rule[0cm]{0cm}{.25cm}
            &&        & Rmse         &$0.0326$&$0.0370$&$0.0288$&$0.0704$&$0.0832$&$0.0055$ &&&\\
            
            &&&&&\\
            &Student& $500$& Mean         &$0.1511$&$-0.1956$&$0.3820$&$0.3332$&$-0.0013$&$0.00585$
 &&\multirow{2}*{$0.070$}&\\
\rule[0cm]{0cm}{.25cm}
            &&        & Rmse        &$0.0445$&$0.0711$&$0.0630$&$0.1934$&$0.0286$&$0.0115$ &&&\\

           \rule[0cm]{0cm}{.4cm}
            && $1000$& Mean         &$0.1475$&$-0.1998$&$0.3902$&$0.3071$&$0.0004$&$0.0049$
            &&\multirow{2}*{$0.055$}&\\
\rule[0cm]{0cm}{.25cm}
            &&        & Rmse         &$0.0302$&$0.0400$&$0.0545$&$0.1394$&$0.0180$&$0.0092$ &&&\\

\hdashline[2pt/3pt]                 

\rule[0cm]{0cm}{.3cm}
  \rule[0cm]{0cm}{.25cm}
      {\bf S}$_1$ &Gaussian& $500$& Mean         &$0.1520$&$-0.2012$&$0.3857$&$0.2997$&$0.0806$&$0.0049$
 &&&\multirow{2}*{$0.750$}\\
\rule[0cm]{0cm}{.25cm}
            &&        & Rmse         &$0.0488$&$0.0547$&$0.0439$&$0.0902$&$0.0316$&$0.0095$ &&&\\

           \rule[0cm]{0cm}{.4cm}
            && $1000$& Mean         &$0.1511$&$-0.1996$&$0.3889$&$0.3006$&$0.0793$&$0.0031$ 
            &&&\multirow{2}*{$0.970$}\\
\rule[0cm]{0cm}{.23cm}
            &&        & Rmse         &$0.0327$&$0.0365$&$0.0295$&$0.0672$&$0.0212$&$0.0055$ &&&\\
   
   &&&&&\\
            &Student& $500$& Mean         &$0.1518$&$-0.2004$&$0.3790$&$0.3149$&$0.0803$&$0.0089$
 &&&\multirow{2}*{$0.825$}\\
\rule[0cm]{0cm}{.25cm}
            &&        & Rmse        &$0.0396$&$0.0616$&$0.0695$&$0.1734$&$0.0275$&$0.0178$ &&&\\

           \rule[0cm]{0cm}{.4cm}
            && $1000$& Mean         &$0.1497$&$-0.1977$&$0.3862$&$0.3076$&$0.0805$&$0.0047$
            &&&\multirow{2}*{$0.990$}\\
\rule[0cm]{0cm}{.25cm}
            &&        & Rmse         &$0.0295$&$0.0452$&$0.0512$&$0.1334$&$0.0198$&$0.0090$ &&&\\
           
 \Xhline{0.6pt}
 
 \rule[0cm]{0cm}{.3cm}
  \rule[0cm]{0cm}{.25cm}
  
    {\bf S'}$_0$&Gaussian& $500$& Mean         &$1.0003$&$0.4009$&$0.4793$&$0.2004$&$0.0020$&$0.0065$
 &&\multirow{2}*{$0.030$}&\\
\rule[0cm]{0cm}{.25cm}
            &&        & Rmse        &$0.0815$&$0.0447$&$0.0672$&$0.0285$&$0.0383$&$0.0126$&&&\\

           \rule[0cm]{0cm}{.5cm}
            && $1000$& Mean         &$1.0015$&$0.4005$&$0.4854$&$0.1986$&$-0.0012$&$0.0050$
             &&\multirow{2}*{$0.055$}&\\
\rule[0cm]{0cm}{.23cm}
            &&        & Rmse         &$0.0481$&$0.0277$&$0.0485$&$0.0200$&$0.0295$&$0.0095$ &&&\\

        &&&&&\\
            &Student& $500$& Mean    &$1.0061$&$0.3957$&$0.4734$&$0.2011$&$0.0005$&$0.0101$
 &&\multirow{2}*{$0.075$}&\\
\rule[0cm]{0cm}{.25cm}
            &&        & Rmse        &$ 0.0830$&$0.0465$&$0.1230$&$0.0583$&$0.0393$&$0.0197$ &&&\\

           \rule[0cm]{0cm}{.4cm}
            && $1000$& Mean         &$1.0025$&$0.3967$&$0.4852$&$0.1976$&$0.0003$&$0.0088$
            &&\multirow{2}*{$0.065$}&\\
\rule[0cm]{0cm}{.25cm}
            &&        & Rmse         &$0.0554$&$0.0296$&$0.0880$&$0.0349$&$0.0245$&$0.0181$ &&&\\
                
       \hdashline[2pt/3pt]
            
\rule[0cm]{0cm}{.3cm}
  \rule[0cm]{0cm}{.25cm}
       {\bf S'}$_1$ &Gaussian& $500$& Mean        &$1.0022$&$0.3993$&$0.4943$&$0.2005$&$0.0704$&$0.0702$
 &&&\multirow{2}*{$0.795$}\\
\rule[0cm]{0cm}{.25cm}
            &&        & Rmse         &$0.0821$&$0.0443$&$0.0803$&$0.0293$&$0.0449$&$0.0315$ &&&\\

           \rule[0cm]{0cm}{.5cm}
            && $1000$& Mean         &$1.0019$&$0.3970$&$0.5052$&$0.1979$&$0.0680$&$0.0672$
 &&&\multirow{2}*{$0.975$}\\
\rule[0cm]{0cm}{.25cm}
            &&        & Rmse         &$0.0587$&$0.0310$&$0.0548$&$0.0217$&$0.0322$&$0.0224$ &&&\\

       &&&&&\\
            &Student& $500$& Mean         &$1.0098$&$0.3979$&$0.5015$&$0.1956$&$0.0711$&$0.0759$
             &&&\multirow{2}*{$0.665$}\\
\rule[0cm]{0cm}{.25cm}
            &&        & Rmse        &$0.0837$&$0.0461$&$0.1575$&$0.0576$&$0.0431$&$0.0522$ &&&\\

           \rule[0cm]{0cm}{.4cm}
            && $1000$& Mean         &$0.9996$&$0.3990$&$0.4949$&$0.2034$&$0.0733$&$0.0692$
             &&&\multirow{2}*{$0.930$}\\\\
\rule[0cm]{0cm}{.25cm}
            &&        & Rmse         &$0.0576$&$0.0318$&$0.1119$&$0.0420$&$0.0295$&$0.0399$&&&\\
            
             \Xhline{0.6pt}
 
 \rule[0cm]{0cm}{.3cm}
  \rule[0cm]{0cm}{.25cm}
    {\bf S''}$_0$    &Gaussian& $500$& Mean         &$37.9776$&$0.3259$&$32.2873$&$0.0198$&$0.0025$&
 &&\multirow{2}*{$0.085$}&\\
\rule[0cm]{0cm}{.25cm}
            &&        & Rmse        &$2.4524$&$0.0373$&$2.8635$&$0.0021$&$0.0350$&&&&\\

           \rule[0cm]{0cm}{.5cm}
            && $1000$& Mean         &$37.7283$&$0.3274$&$31.9529$&$0.0200$&$0.0042$&
             &&\multirow{2}*{$0.060$}&\\
\rule[0cm]{0cm}{.23cm}
            &&        & Rmse         &$2.2560$&$0.0302$&$2.7856$&$0.0018$&$0.0325$&&&&\\
            
       &&&&&\\
            &Student& $500$& Mean         &$38.0380$&$0.3254$&$32.0056$&$0.0194$&$0.0027$
 &&&\multirow{2}*{$0.065$}&\\
\rule[0cm]{0cm}{.25cm}
            &&        & Rmse        &$2.3269$&$0.04313$&$2.8796$&$0.0036$&$0.0370$ &&&\\

           \rule[0cm]{0cm}{.4cm}
            && $1000$& Mean         &$37.8905$&$0.3279$&$31.7781$&$0.0199$&$0.0019$
            &&&\multirow{2}*{$0.040$}&\\
\rule[0cm]{0cm}{.25cm}
            &&        & Rmse         &$2.1403$&$0.0273$&$2.8731$&$0.0026$&$0.0291$&&&&\\     
            
       \hdashline[2pt/3pt]
            
\rule[0cm]{0cm}{.3cm}
  \rule[0cm]{0cm}{.25cm}
       {\bf S''}$_1$&Gaussian& $500$& Mean         &$37.9698$&$0.3265$&$33.0176$&$0.0189$&$-0.2092$&
 &&&\multirow{2}*{$0.985$}\\
\rule[0cm]{0cm}{.25cm}
            &&        & Rmse         &$2.4005$&$0.0401$&$2.6480$&$0.0040$&$0.0310$& &&&\\

           \rule[0cm]{0cm}{.5cm}
            && $1000$& Mean         &$38.1928$&$0.3255$&$32.7826$&$0.0193$&$-0.2113$&
 &&&\multirow{2}*{$1.000$}\\
\rule[0cm]{0cm}{.25cm}
            &&        & Rmse         &$2.1828$&$0.0289$&$2.5628$&$0.0031$&$0.0263$&&&&\\
            
      &&&&&\\
            &Student& $500$& Mean         &$38.2794$&$0.3232$&$32.5255$&$0.0194$&$-0.2117$
 &&&&\multirow{2}*{$0.995$}\\
\rule[0cm]{0cm}{.25cm}
            &&        & Rmse            &$2.2897$&$0.0397$&$2.8555$&$0.0065$&$0.0282$&&&&\\

           \rule[0cm]{0cm}{.4cm}
            && $1000$& Mean         &$38.3885$&$0.3288$&$32.1560$&$0.0198$&$-0.2147$
            &&&&\multirow{2}*{$1.000$}\\
\rule[0cm]{0cm}{.25cm}
            &&        & Rmse         &$2.1398$&$0.0258$&$2.7629$&$0.0046$&$0.02653$&&&&\\

 \Xhline{.7pt}
\end{tabular}
\end{table}

\section{Real data example}\label{sect_Real_data}
 We consider the daily concentrations of the PM$_{10}$ (particulate  matter  with a diameter less than 10$\mu m$) in the Vit\'oria metropolitan area; see Figure \ref{Graphe_PM10_Residuals}(a).
    These data as well as those of the other meteorological variables (see below) are obtained from the State Environment and Water Resources Institute, and were collected at eight monitoring stations.
 We focus on the data from January 21st, 2005 to March 04th 2006, these observations on 408 days are a part of a large dataset (available at https://rss.onlinelibrary.wiley.com/pb-assets/hub-assets/rss/Datasets/RSSC\%2067.2/C1239deSouza-1531120585220.zip) which were analyzed by \nocite{Souza2018} Souza {\it et al.} (2018) to quantify the association between respiratory disease and air pollution concentrations.

\medskip

The variables considered are: the average concentration for the particulate matter (PM$_{10}$, $\mu g m^{-3}$), sulphur dioxide (SO$_2$, $\mu g m^{-3}$), nitrogen dioxide (NO$_2$,  $\mu g m^{-3}$), carbon monoxide (CO, $\mu g m^{-3}$), ozone (O$_3$, $\mu g m^{-3}$); and  Air relative humidity (RH, \%). 
Table \ref{Table_element_stat} displays some elementary statistics of these variables.

\begin{table}[h!]
\scriptsize
\centering
\caption{ \it \small Some elementary statistics of the variables PM$_{10}$, SO$_2$, NO$_2$, CO, O$_3$, RH, for the period from January 21st, 2005 to March 04th 2006. }
\label{Table_element_stat}
\vspace{.1cm}
\begin{tabular}{cccccccc}
\hline 
 Variable &Mean&SD&Min&$Q_1$&Med&$Q_3$&Max\\
\Xhline{.6pt}
\rule[0cm]{0cm}{.35cm}
 PM$_{10}$ ($\mu g m^{-3}$)   & 32.04 & 8.62 & 11.16 & 26.19 & 31.87 & 36.92 & 66.60 \\
           &  &  &  &  &  &  &   \\       
  SO$_2$ ($\mu g m^{-3}$)  & 11.64 & 2.74 & 4.89 & 9.75 & 11.63 & 13.48 &19.29\\ 
             &  &  &  &  &  &  &   \\       
    NO$_2$ ($\mu g m^{-3}$)  & 24.57 & 6.41 & 10.47 & 19.93 & 23.79 & 28.84 & 46.84\\ 
            &  &  &  &  &  &  &   \\       
    CO ($\mu g m^{-3}$)    & 969.70 & 254.85 & 456.00 & 785.10 & 951.00 & 1129.70 &2141.50\\
            &  &  &  &  &  &  &   \\       
      O$_3$ ($\mu g m^{-3}$)  & 29.96   & 7.68 & 16.76 & 24.79 & 28.46 & 33.96 & 66.52 \\
            &  &  &  &  &  &  &   \\       
      RH (\%)  & 79.21  & 6.52 & 62.45 & 74.36 & 78.67 & 83.62 & 95.39 \\ 
 \Xhline{.9pt}
\end{tabular}
\end{table} 

\medskip

As pointed out by Ng and Awang (2018), PM$_{10}$ is a notorious air pollutant associated in particular with detrimental health impacts; it affects the respiratory and cardiopulmonary functions and increases the morbidity and mortality rate of related diseases.
Therefore, forecasting the PM$_{10}$ concentration and understanding its relation with other factors is an important issue.
Several models, including among others models, ARIMA, MLR (multiple linear regression), RTSE (Regression with time series error), were considered; we refer to Ng and Awang (2018) and  Ng (2017) and the references therein for an overview of this issue.

\medskip

In this section, we focus on the forecast of the PM$_{10}$ concentration from some meteorological variables of the previous day.  We apply the model (\ref{FDAR_model}) with $p_1, p_2, q_1, q_2 \in \{ 0 , 1\}$, and the covariate $X_t = (SO_{2,t}, NO_{2,t}, CO_t, O_{3,t}, RH_t)'$ (the value of the corresponding variable at day $t$). The following issues are addressed.
\begin{enumerate}
\item \textbf{Model selection}.  The aim is to select the orders $p_1, p_2, q_1, q_2$ and the "best" subset of the covariates that are the major factor related to the next-day PM$_{10}$ concentration. For this purpose, we consider all the combination of the covariates with $p_1, p_2, q_1, q_2 \in \{ 0 , 1\}$; which represents a collection of 376 models.
The procedure based on the penalized criteria $\widehat{C}_n(m)$ (see (\ref{Cont_pen})) is applied with the regularization parameter $\kappa_n = \log n, ~ 2\log \log n, ~ 3.5 \log \log n, ~ 5 \log \log n$. These criteria (BIC and HQC) select the model with $p_1= p_2=q_1=1$, $q_2=0$ and the covariate $RH$. Compared to ARIMA, this model is also preferred.
This result is in accordance with some existing works (see, for instance, Ng and Awang (2018) and  Ng (2017)) which have found that the air humidity of the previous day is an important factor related to PM$_{10}$ concentration.   
\item \textbf{Estimation and significance test}. The estimated model is:
\[ PM_{10,t} = \underset{(2.671)}{37.946} +   \underset{(0.024)}{0.330} PM_{10,t-1}  - \underset{(0.028)}{0.210} RH_{t-1}+\xi_t \sqrt{\underset{(3.362)}{32.108} + \underset{(0.003)}{0.023} PM_{10,t-1}^2 }, \]
where in parentheses are the standard errors of the estimators obtained from the robust sandwich matrix.
The test (\ref{Test}) with $\vartheta_0=0$ is now applied for testing the significance of the covariate. At the nominal level $\alpha=0.05$, the critical value of the test, computed from $(\Gamma Z^{\mathcal C})'(\Gamma \Sigma \Gamma')^{-1}\Gamma Z^{\mathcal C}$ is 2.68 and the statistic, computed from $W_n$ (see (\ref{Stat_Wn})) is 10.86. Thus, the null hypothesis is rejected. 
Figure \ref{Graphe_PM10_Residuals} displays the histogram and cumulative periodogram of the residuals as well as the autocorrelation functions of the squared residuals. From these findings, the residuals do not show any signs of correlation.

\begin{figure}[h!]
\begin{center}
\includegraphics[height=9cm, width=15cm]{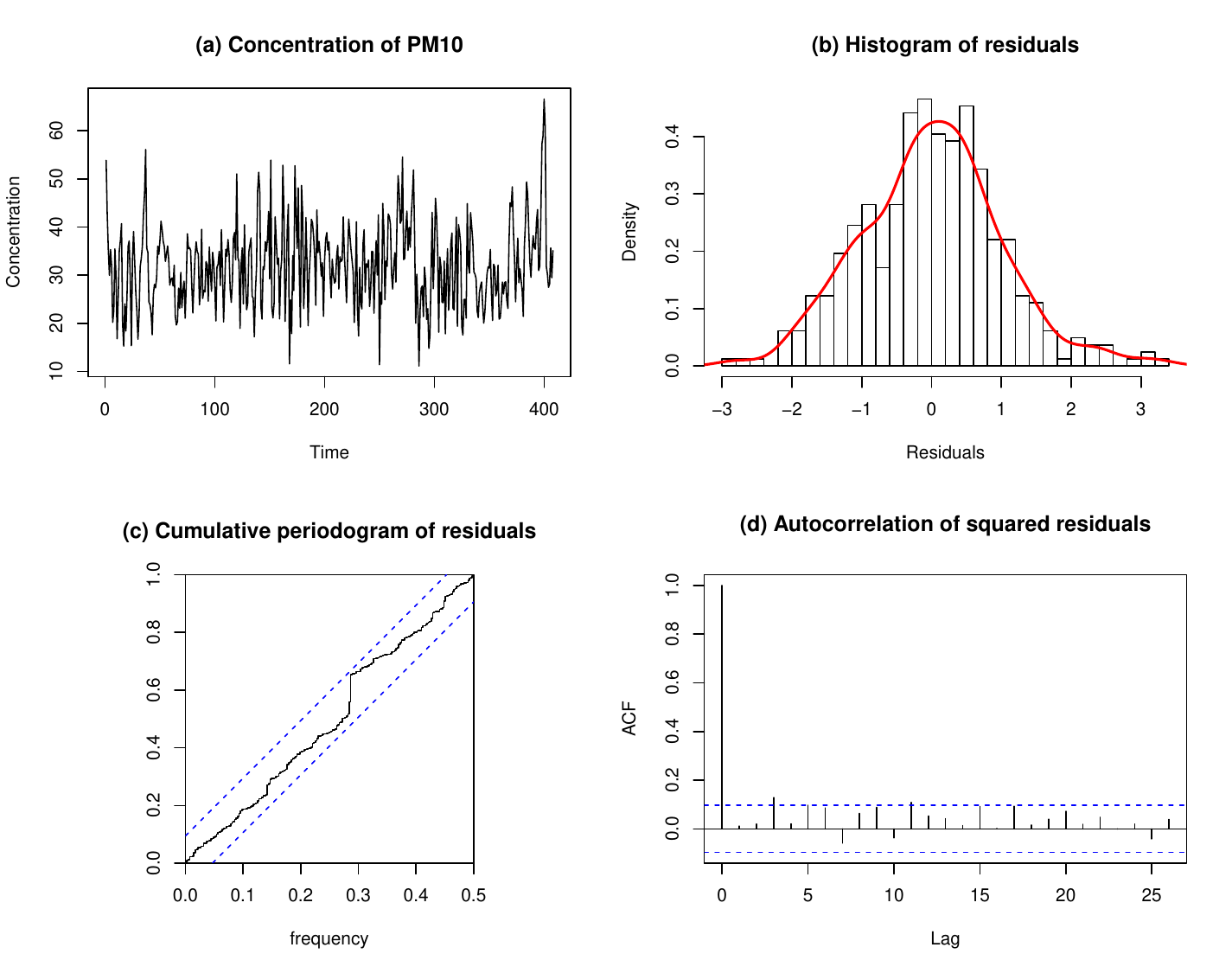}
\end{center}
\vspace{-.6cm}
\caption{\it (a) The daily concentrations of the PM$_{10}$ from January 21st, 2005 to March 04th 2006; (b), (c) the histogram and the cumulative periodogram of the residuals; (d) the autocorrelation functions of the squared residuals.}
\label{Graphe_PM10_Residuals}
\end{figure}
\end{enumerate}
In conclusion of this section, let us stress that, several authors have applied the MLR and the RTSE on other data and have found that NO$_2$, CO, O$_3$    were also important factors associated to the PM$_{10}$ concentration (see Ng and Awang (2018) and  Ng (2017) and the references therein). For the data considered here, and by applying the model (\ref{FDAR_model}), it appears that these variables were less important in for forecasting the PM$_{10}$ concentrations.

\section{Summary and conclusion}
This paper considers a general class of causal processes with exogenous covariates in a semiparametric framework. This class is quite extensive and many classical processes such as ARMA-GARCH, ARMAX-GARCH-X, APARCH-X,$\cdots$ are particular cases. 
Sufficient conditions for the existence of a stationary and ergodic solution are provided. \\
\noindent A quasi likelihood estimator is performed for inference; the consistency of this estimator is established and the asymptotic distribution is derived. This distribution coincides with the Gaussian one, when the true parameter is an interior point of the parameter's space.\\
\noindent A Wald-type statistic is proposed for testing the significance test of parameter. The asymptotic studies show that, this test has correct size asymptotically and is consistent in power. In certain cases, this test can be used in particular to test the relevance of the exogenous covariates. \\
\noindent The model selection question for the class $\mathcal {AC}$-$X(M_{\theta^*},f_{\theta^*})$ is carried out by a penalized quasi likelihood contrast. The weak and the strong consistency of the proposed procedure is established. These results provide sufficient conditions for the consistency of the BIC and the HQC procedures. Simulation study shows that, the empirical and the theoretical results are overall in accordance.\\
\indent An extension of this works is to address the inference, the significance test of parameter, the  model selection problem for the class $\mathcal {AC}$-$X(M_{\theta^*},f_{\theta^*})$ with a non Gaussian quasi likelihood. For instance, as pointed out by Kengne (2021), the use of the Laplacian quasi likelihood will allow to reduce the order of moments imposed to the process.
Other topics of a research project, are the change-point detection and the prediction question (see for instance Ing (2003), Ing and Wei (2003, 2005))  for this class of models.


\section{Proofs of the main results}
 To simplify the expressions, in the proofs of Theorems \ref{th1}, \ref{th2} and \ref{th4}, we will use the conditional Gaussian quasi-log-likelihood given by $L_n(\theta)=  -\sum\limits_{t=1}^{n}  q_t(\theta)$ and $\widehat  L_n(\theta)=  -\sum\limits_{t=1}^{n} \widehat q_t(\theta)$. 
Throughout the sequel, $C$ denotes a positive constant whom value may differ from an inequality to another.
 
 \subsection{Proof of Proposition \ref{prop1}}
 We verify that the process $Z_t := (Y_t; X_t)$ satisfies the conditions required for the Theorem 3.1 in Doukhan and Wintenberger \cite{Doukhan2008}.   
 According to (\ref{Model}), for all $t \in \Z$,
 \begin{align*}
 Z_t 
 &=
  \big(M_{\theta^*}(Y_{t-1}, \ldots; X_{t-1},\ldots)\xi_t + 
      f_{\theta^*}(Y_{t-1}, \ldots; X_{t-1},\ldots) ; \,g(X_{t-1},\ldots; \eta_t)\big)\\
 &=
 F(Z_{t-1},Z_{t-2},\ldots;U_t),
 \end{align*}
 with $U_t=(\xi_t,\eta_t )$ and $F(\pmb{z};U_t) = \big(M_{\theta^*}(y_1, \ldots; x_1,\ldots)\xi_t + 
      f_{\theta^*}(y_1, \ldots; x_1,\ldots); \,g(x_1,\ldots; \eta_t) \big)$
for all  $\pmb{z} = \big((y_k,x_k) \big)_{k \in \N} \in ( \R^{d_x+1})^\N$.
 Thus, the equation (1.1) of \cite{Doukhan2008} holds for $(Z_t)_{t \in \Z}$. 
 For a vector $z=(y,x) \in \R^{d_x+1}$, define the norm $\| z \|_w = |y| + w_x \| x \| $ for some $ w_x > 0$. 
 According to Doukhan and Wintenberger (2008), it suffices to show that:
 \begin{itemize}
  \item[(i)] $\E \| F(\pmb{z};U_0) \|^r_w < \infty $ for some $\pmb{z} \in ( \R^{d_x+1})^\N$;
  \item[(ii)] there exists a non-negative sequence $(\alpha_k(F))_{k \geq 1}$ satisfying $\sum_{k \geq 1} \alpha_k(F) < 1$ such that, for all $\pmb{z}, \tilde{\pmb{z}} \in (\R^{d_x+1})^\N$,
  \[\E \| F(\pmb{z};U_0) - F(\tilde{\pmb{z}};U_0) \|^r_w \leq \sum_{k \geq 1} \alpha_k(F) \|z_k - \tilde{z}_k \|_w  . \]
\end{itemize}  
 
  \medskip
    
    \noindent
  Using the condition (\ref{exp_Lip_cov}), the part (i) is directly obtained from the assumptions \textbf{A}$_0(f_\theta,\Theta)$, \textbf{A}$_0(M_\theta,\Theta)$.\\
 To prove (ii), let $\pmb{z}=(z_1,\ldots),\ \tilde{\pmb{z}}=(\tilde{z}_1,\ldots) \in ( \R^{d_x+1})^\N$  such that $z_k=(y_k,x_k)$ and $\tilde{z}_k=(\tilde{y}_k,\tilde{x}_k)$ for all $k\geq 1$.
From \textbf{A}$_0(f_\theta,\Theta)$, \textbf{A}$_0(M_\theta,\Theta)$ and (\ref{exp_Lip_cov}), we get 
  \begin{align*}
  &\left\|\left\|F(\pmb{z};U_0)-F(\tilde{\pmb{z}};U_0)\right\|_w\right\|_r \\
  &\hspace{.2cm}\leq
  \big\|\left\|\left(M_\theta(y_{1}, \ldots; x_{1},\ldots) - M_\theta(\tilde{y}_{1}, \ldots; \tilde{x}_{1},\ldots)\right)\xi_0\right\|_\Theta 
  +
  \left\|f_\theta(y_{1}, \ldots; x_{1},\ldots) - f_\theta(\tilde{y}_{1}, \ldots; \tilde{x}_{1},\ldots)\right\|_\Theta \big\|_r\\
  &\hspace{1.5cm} +
  w_x \left\|g(x_{1},\ldots; \eta_0)-g(\tilde{x}_{1},\ldots; \eta_0) \right\|_r\\
  &\hspace{.2cm}\leq
    \sum\limits_{k=1}^{\infty} \left(\alpha^{(0)}_{k,Y}(f_\theta,\Theta) + \|\xi_0\|_r  \alpha^{(0)}_{k,Y}(M_\theta,\Theta)\right) |y_k-\tilde{y}_k|+ \sum\limits_{k=1}^{\infty} \left(\alpha^{(0)}_{k,X}(f_\theta,\Theta) + \|\xi_0\|_r  \alpha^{(0)}_{k,X}(M_\theta,\Theta)\right)\|x_k-\tilde{x}_k\|\\
  &\hspace{1.5cm} +
w_x \sum\limits_{k=1}^{\infty} \alpha_k(g) \left\| x_k-\tilde{x}_k \right\|\\  
 &\hspace{.2cm}\leq
 \sum\limits_{k=1}^{\infty} \left(\alpha^{(0)}_{k,Y}(f_\theta,\Theta) + \|\xi_0\|_r   \alpha^{(0)}_{k,Y}(M_\theta,\Theta)\right) |y_k-\tilde{y}_k|\\
 &\hspace{1.5cm} 
 + w_x \sum\limits_{k=1}^{\infty} \Big(\frac{1}{w_x} \big( \alpha^{(0)}_{k,X}(f_\theta,\Theta) + \|\xi_0\|_r   \alpha^{(0)}_{k,X}(M_\theta,\Theta)\big)+\alpha_k(g)\Big)\|x_k-\tilde{x}_k\|\\
  &\hspace{.2cm} \leq
  \sum\limits_{k=1}^{\infty} \alpha_k(F) \|z_k-\tilde{z}_k\|_{w}
 \end{align*}
 with 
 $
 \alpha_k(F) = \max \left\{ \alpha^{(0)}_{k,Y}(f_\theta,\Theta) + \|\xi_0\|_r \alpha^{(0)}_{k,Y}(M_\theta,\Theta),\ \frac{1}{w_x} \big( \alpha^{(0)}_{k,X}(f_\theta,\Theta) + \|\xi_0\|_r \alpha^{(0)}_{k,X}(M_\theta,\Theta)\big)+\alpha_k(g)\right\}.
 $
 Thus, to get $\sum\limits_{k=1}^{\infty} \alpha_k(F)<1$,  it suffices to  choose $w_x$ sufficiently large, such that 
 \[
  w_x > \frac{\sum_{k\geq 1} \left\{\alpha^{(0)}_{k,X}(f_\theta,\Theta) + \|\xi_0\|_r \alpha^{(0)}_{k,X}(M_\theta,\Theta)\right\}}{1- \sum_{k\geq 1} \max\left\{ \alpha_k(g),\,\alpha^{(0)}_{k,Y}(f_\theta,\Theta) + \|\xi_0\|_r \alpha^{(0)}_{k,Y}(M_\theta,\Theta)\right\}}.
 \]
 This completes the proof of the proposition.
     \begin{flushright}
$\blacksquare$
 \end{flushright}
 
  
 \subsection{Proof of Theorem \ref{th1}}\label{proof_th1}
 We consider the following lemma. 
 
 \begin{lem}\label{lem1} 
 Assume that the assumptions of Theorem \ref{th1} hold. Then
 
\begin{equation*}
\frac{1}{n} \big\|\widehat L_n(\theta) - L_n(\theta)  \big\|_\Theta \limitepsn 0 .
 \end{equation*}
 \end{lem}

 \noindent
 \emph{\bf Proof of Lemma \ref{lem1}}\\
Remark that 
\begin{align*}
  \frac{1}{n} \| \widehat{L}_n(\theta)-L_n(\theta)\|_\Theta
 &\leq
 \frac{1}{n} \sum_{t=1}^{n}\|\widehat{q}_{t}(\theta)-q_t(\theta) \|_\Theta.
 \end{align*}
Hence, by Corollary 1 of  Kounias and Weng (1969), with $2\leq \tilde  r \leq \min\{3,r\}$ (without loss of generality), it suffices to show that 
\begin{equation}\label{Kounias_lem1}
 \sum_{\ell\geq 1}^{}\frac{1}{\ell^{\tilde r/3}}\E \left(\|\widehat{q}_{\ell}(\theta)-q_\ell(\theta) \|^{\tilde r/3}_\Theta\right) < \infty.
\end{equation}
For any $\theta \in \Theta$, by applying the mean value theorem at the functions  $x \mapsto \frac{1}{x^2}$ and  $x \mapsto \log x$, we have
\begin{align*}
|\widehat{q}_{t}(\theta)-q_t(\theta)| 
&\leq
\Big|\frac{(Y_t-\widehat f^t_\theta)^2}{\widehat H^t_\theta}-\frac{ (Y_t-f^t_\theta)^2}{H^t_\theta}\Big|+|\log \widehat H^t_\theta-\log H^t_\theta| \\
&\leq
\Big|\frac{(Y_t-\widehat f^t_\theta)^2}{ (\widehat M^t_\theta)^2}-\frac{ (Y_t-f^t_\theta)^2}{(M^t_\theta)^2}\Big|+2\big|\log |\widehat M^t_\theta|-\log |M^t_\theta|\big| \\
&\leq
\Big|(Y_t-\widehat f^t_\theta)^2 \Big(\frac{1}{(\widehat M^t_\theta)^2}-\frac{1}{(M^t_\theta)^2}\Big)+\frac{1}{(M^t_\theta)^2} \left( (Y_t-\widehat f^t_\theta)^2-(Y_t- f^t_\theta)^2\right)\Big|+ \frac{2}{\underline{h}^{1/2}}| \widehat M^t_\theta-  M^t_\theta| \\
&\leq
\frac{2}{\underline{h}^{3/2}} |Y_t-\widehat f^t_\theta|^2  |\widehat M^t_\theta-M^t_\theta|
+\frac{1}{\underline{h}} |\widehat f^t_\theta-f^t_\theta| |\widehat f^t_\theta+ f^t_\theta -2Y_t|+ \frac{2}{\underline{h}^{1/2}}| \widehat M^t_\theta-  M^t_\theta| \\
&\leq
C \big( (|Y_t-f^t_\theta|^2 +1) |\widehat M^t_\theta- M^t_\theta| +
 |\widehat f^t_\theta-f^t_\theta| |\widehat f^t_\theta+f^t_\theta-2Y_t|\big).
\end{align*}
This implies
\begin{multline*}
\E\big[\|\widehat q_t(\theta) - q_t(\theta)\|^{\tilde r/3}_\Theta\big]
\leq
C \Big( \E \Big[\big(\|Y_t-f^t_\theta\|^2_\Theta+1\big)^{\tilde r/3} \|\widehat M^t_\theta- M^t_\theta\|^{\tilde r/3}_\Theta\Big]\\
+
\E \Big[\|\widehat f^t_\theta-f^t_\theta\|^{\tilde r/3}_\Theta \big(\|\widehat f^t_\theta\|_\Theta+\|f^t_\theta\|_\Theta+2|Y_t|\big)^{\tilde r/3}\Big]
\Big).
\end{multline*}
Moreover, since $ \theta^* \in \Theta(r)$ for some $r \geq 2$, by Assumption \textbf{A}$_0(\Psi_\theta,\Theta)$, one can easily show that:
%
\begin{align}
& \label{Res1_A0}
\bullet~~
\E \left[|Y_t|^{r} + \|f^{t}_\theta\|^{r}_\Theta +\|\widehat f^{t}_\theta\|^{r}_\Theta + \|M^{t}_\theta\|^{r}_\Theta + \|\widehat M^{t}_\theta\|^{r}_\Theta + \|H^{t}_\theta\|^{r/2}_\Theta + \|\widehat H^{t}_\theta\|^{r/2}_\Theta\right]< \infty;\\
& \label{Res2_A0}
\bullet~~
     \left\{
\begin{array}{l} 
     \E(\|\widehat f^{t}_\theta-f^{t}_\theta\|^r_\Theta) 
       \leq  C\Big(\sum\limits_{k \geq t}  \big \{ \alpha^{(0)}_{k,Y}(f_\theta,\Theta) +\alpha^{(0)}_{k,X}(f_\theta,\Theta)\big\} \Big)^{r};\\
       \rule[0cm]{0cm}{.6cm}
     \E(\|\widehat M^{t}_\theta-M^{t}_\theta\|^r_\Theta) 
       \leq  C\Big(\sum\limits_{k \geq t}  \big\{ \alpha^{(0)}_{k,Y}(M_\theta,\Theta) +\alpha^{(0)}_{k,X}(M_\theta,\Theta)\big\}\Big)^{r}.
\end{array} 
\right.
\end{align}  
       Then, by the H\"{o}lder's inequality, we have 
      \begin{align*}
      & \E \Big[\|\widehat f^\ell_\theta-f^\ell_\theta\|^{\tilde r/3}_\Theta \big(\|\widehat f^\ell_\theta\|_\Theta+\|f^\ell_\theta\|_\Theta+2|Y_\ell|\big)^{\tilde r/3}\Big]\\
      &\hspace{4cm}\leq 
      \Big(\E \big[\|\widehat f^\ell_\theta-f^\ell_\theta\|^{\tilde r}_\Theta\big]\Big)^{1/3} 
      \Big(\E \big[\|\widehat f^\ell_\theta\|_\Theta+\|f^\ell_\theta\|_\Theta+2|Y_\ell|\big]^{\tilde r/2}\Big)^{2/3} \\
      &\hspace{4cm}\leq  
      C  \Big(\sum\limits_{k \geq t}  \big\{ \alpha^{(0)}_{k,Y}(f_\theta,\Theta) +\alpha^{(0)}_{k,X}(f_\theta,\Theta)\big\} \Big)^{\tilde r/3}.
      \end{align*}
 Again, by the H\"{o}lder's inequality, from (\ref{Res1_A0}) and  (\ref{Res2_A0}), we obtain 
 \begin{align*}    
 \E \Big[\big(\|Y_t-f^\ell_\theta\|^2_\Theta+1\big)^{\tilde r/3} \|\widehat M^\ell_\theta- M^\ell_\theta\|^{\tilde r/3}_\Theta\Big]
 &\leq
 \E \Big[\|Y_t+f^\ell_\theta+1\|^{2 \tilde r/3}_\Theta \|\widehat M^\ell_\theta- M^\ell_\theta\|^{\tilde r/3}_\Theta\Big]\\
 &\leq 
  \big(\E \big[\|Y_t+f^\ell_\theta+1\|^{\tilde r}_\Theta\big]\big)^{2/3}  \big(\E \big[\|\widehat M^\ell_\theta- M^\ell_\theta\|^{\tilde r}_\Theta\big]\big)^{1/3}\\
 &\leq 
 C \Big(\sum\limits_{k \geq t}  \big\{ \alpha^{(0)}_{k,Y}(M_\theta,\Theta) +\alpha^{(0)}_{k,X}(M_\theta,\Theta)\big\} \Big)^{\tilde r/3}.
 \end{align*}
Hence, from (\ref{cond_th1}), we deduce
 \begin{align*}
&\sum_{\ell\geq 1}^{}\frac{1}{\ell^{\tilde r/3}}\E \left(\|\widehat{q}_{\ell}(\theta)-q_\ell(\theta) \|^{\tilde r/3}_\Theta\right)\\
&\hspace{2cm}\leq
C\sum_{\ell\geq 1}^{}\frac{1}{\ell^{\tilde r/3}} \Big( \sum\limits_{k \geq \ell}  \big\{ \alpha^{(0)}_{k,Y}(f_\theta,\Theta) +\alpha^{(0)}_{k,X}(f_\theta,\Theta)+
 \alpha^{(0)}_{k,Y}(M_\theta,\Theta) +\alpha^{(0)}_{k,X}(M_\theta,\Theta)\big\} \Big)^{\tilde r/3}\\
&\hspace{2cm} \leq  C\sum_{\ell\geq 1}^{}\frac{1}{\ell^{\tilde r/3}} \Big(\dfrac{1}{\ell^{\gamma-1}} \Big)^{\tilde r/3}\leq   
      C  \sum_{\ell\geq 1}^{}\frac{1}{\ell^{\tilde r \gamma/3}} <\infty,
 \end{align*}
  where the last inequality holds since  $\gamma> 3/2$. 
  Thus, the condition (\ref{Kounias_lem1}) is satisfied. This completes the proof Lemma \ref{lem1}.
  \begin{flushright}
 $~~~~~~\blacksquare$ \\
\end{flushright}

\noindent 
To complete the proof of Theorem \ref{th1}, we will show that: (1.) $\E \left[\left\| q_{t}(\theta)\right\|_\Theta\right] <\infty$ and (2.) the function $\theta \mapsto -\E[q_{0}(\theta)]$ has a unique maximum at $\theta^*$.  
\begin{enumerate}
	\item [(1.)] For all $\theta \in \Theta$, using the inequality $|\log(x)| \leq  |x - 1|$ for all $x >1$, we have
	\begin{align*} 
	| q_{t}(\theta)|  
	&\leq
	\frac{1}{H^t_\theta}|Y_t-f^t_\theta)|^2+\big|\log \big(\frac{H^t_\theta}{\underline{h}}\big) + \log(\underline{h})\big|\\
		&\leq
	\frac{1}{\underline{h}}(Y^2_t+(f^t_\theta)^2 +2 Y_tf^t_\theta)+\big|\frac{H^t_\theta}{\underline{h}}-1\big| + |\log(\underline{h})|\\
		&\leq
	C\left(Y^2_t+(f^t_\theta)^2 +2 Y_tf^t_\theta+|M^t_\theta|^2\right)+ C.
	\end{align*}
	Hence, from (\ref{Res1_A0}), we deduce 
	\[
	\E \left[\left\| q_{t}\right\|_\Theta\right]
	\leq
	C \Big( \E [Y^2_t]+\E \|f^t_\theta\|^2_\Theta +2 \big( (\E [Y^2_t])^{1/2}(\E [\|f^t_\theta\|^2_\Theta] )^{1/2}\big) +\E [\|M^t_\theta\|^2_\Theta] \Big)+C <	\infty,
	\]
	which shows that (1.) holds.
	
		\item [(2.)]
			Let $\theta \in \Theta$ with $\theta \neq \theta^*$. We have 
			\begin{equation}\label{eq1_proof_th1}
			\E[q_{0}(\theta)]-\E[q_{0}(\theta^*)]=\E \big[\E\left[\left(q_{0}(\theta)-q_{0}(\theta^*) \right)| \mathcal{F}_{-1}\right]\big].
			\end{equation}
			Moreover,
			\begin{align*}
			\E\left[\left(q_{0}(\theta)-q_{0}(\theta^*) \right)| \mathcal{F}_{-1}\right]
			&=
			\E \big[ \frac{(Y_0-  f^0_\theta)^2}{H^0_\theta}+\log H^0_\theta-
    \frac{(Y_0-  f^0_{\theta^*})^2}{H^0_{\theta^*}}-\log H^0_{\theta^*}	\big| \mathcal{F}_{-1} \big]\\
     &=
    -\log \big( \frac{H^0_{\theta^*}}{H^0_{\theta}}\big)  +  \frac{\E \left[(Y_0-  f^0_\theta)^2 \big| \mathcal{F}_{-1} \right]}{  H^0_\theta}-\frac{\E \left[(Y_0-  f^0_{\theta^*})^2 \big| \mathcal{F}_{-1} \right]}{  H^0_{\theta^*}}\\
    &=
    -\log \big(\frac{H^0_{\theta^*}}{H^0_{\theta}}\big)  -1+  \frac{\E \left[(Y_0-  f^0_{\theta^*}+  f^0_{\theta^*} - f^0_\theta)^2 \big| \mathcal{F}_{-1} \right]}{H^0_\theta}\\
      &=
    \frac{H^0_{\theta^*}}{H^0_{\theta}}-\log \big(\frac{H^0_{\theta^*}}{H^0_{\theta}}\big)  -1+  \frac{( f^0_{\theta^*} - f^0_\theta)^2}{H^0_\theta}.
			\end{align*}
Therefore, using (\ref{eq1_proof_th1}) and by applying the Jensen's inequality, we get
\begin{align*}
\E[q_{0}(\theta)]-\E[q_{0}(\theta^*)]
&=
\E \big[ \frac{H^0_{\theta^*}}{H^0_{\theta}}-\log \big(\frac{H^0_{\theta^*}}{H^0_{\theta}}\big)  -1+  \frac{( f^0_{\theta^*} - f^0_\theta)^2}{H^0_\theta}\big]\\
&\geq
\E \big[ \frac{H^0_{\theta^*}}{H^0_{\theta}}\big]-\log \big(\E \big[\frac{H^0_{\theta^*}}{H^0_{\theta}}\big]\big)  -1+  \E \big[\frac{( f^0_{\theta^*} - f^0_\theta)^2}{H^0_\theta}\big].
\end{align*}
Since $x-\log(x) - 1 > 0$ for any $x > 0,\ x \neq 1$; and $x-\log(x) - 1 = 0$ for $x = 1$, we deduce: 
\begin{itemize}
	\item if $f^0_{\theta^*} \neq f^0_{\theta}~ a.s.$, then $ \E \big[\frac{( f^0_{\theta^*} - f^0_\theta)^2}{H^0_\theta}\big]>0$ and $\E[q_{0}(\theta)]-\E[q_{0}(\theta^*)]>0$,
	
	\item if $f^0_{\theta^*} = f^0_{\theta}~a.s.$, then
	\[
	\E[q_{0}(\theta)]-\E[q_{0}(\theta^*)]=
	\E \big[ \frac{H^0_{\theta^*}}{H^0_{\theta}}-\log \big(\frac{H^0_{\theta^*}}{H^0_{\theta}}\big)  -1\big]
	\]
	 From the identifiability condition (\textbf{A0}), when $\theta^* \neq \theta$ and $f^0_{\theta^*} = f^0_{\theta}~a.s.$, we necessarily have $H^0_{\theta^*} \neq H^0_{\theta}~a.s.$. This implies $\frac{H^0_{\theta^*}}{H^0_{\theta}} \neq 1 ~a.s.$, and thus $\E[q_{0}(\theta)]-\E[q_{0}(\theta^*)]>0$.
\end{itemize}
The equality $\E[q_{0}(\theta)]=\E[q_{0}(\theta^*)]$ holds a.s. if and only if $\theta^* = \theta$. This achieves the proof of (2.).  
	\end{enumerate}
Since $\{(Y_{t},X_t),\,t\in \Z\}$ is stationary and ergodic, the process $\{ q_t(\theta),~t\in \Z\}$ is also a stationary and ergodic sequence. Then, according to (1.), by the uniform strong law of large number applied on the process $\{ q_t(\theta),~ t\in \Z\}$, it holds that 
\[
\big\| \frac{1}{n}L_{n}(\theta) + \E(q_0(\theta)) \big\|_\Theta=  \big\| \frac{1}{n} \sum_{t=1}^{n} q_t(\theta) - \E(q_0(\theta)) \big\|_\Theta \limitepsn 0 .
\]
Then, by Lemma \ref{lem1}, we obtain 
\begin{align} \label{eq2_proof_th1}
\big\|  \frac{1}{n} \widehat L_{n}(\theta) + \E(q_0(\theta)) \big\|_\Theta 
&\leq 
\frac{1}{n}\big\|   \widehat L_{n}(\theta) -  L_{n}(\theta)  \big\|_\Theta + 
\big\|  \frac{1}{n} L_{n}(\theta) + \E(q_0(\theta)) \big\|_\Theta \limitepsn 0.
\end{align}
The part (2.) and (\ref{eq2_proof_th1}) lead to conclude the proof of the theorem.
\begin{flushright}
$\blacksquare$
\end{flushright}


\subsection{Proof of Theorem \ref{th2}}\label{proof_th2}
The following lemma is needed.
\begin{lem}\label{lem2} 
 Assume that the conditions of Theorem \ref{th2} hold. Then
 \rm
 
\begin{enumerate}
	\item [(i.)] 
	$\E\big[ \frac{1}{\sqrt{n}} \big\|\frac{ \partial \widehat L_n(\theta)}{\partial \theta} - \frac{ \partial  L_n(\theta)}{\partial \theta}  \big\|_\Theta \big] \limiten 0$;
	
	\item [(ii.)] 
	$\frac{1}{n} \big\|\frac{ \partial^2 \widehat L_n(\theta)}{\partial \theta \partial \theta'} - \frac{ \partial^2  L_n(\theta)}{\partial \theta \partial \theta'}  \big\|_\Theta  \limitepsn 0$;
	
	\item [(iii.)] $
\big\|\frac{1}{n} \sum\limits_{t=1}^{n}\frac{ \partial^2 q_{t}(\theta)}{\partial \theta \partial \theta'} - \E\big( \frac{ \partial^2   q_{0}(\theta)}{\partial \theta \partial \theta'}\big)  \big\|_\Theta \limitepsn 0 .
$
\end{enumerate}
  \end{lem}

\noindent
\emph{\bf Proof of Lemma \ref{lem2}}
\begin{enumerate}
	\item [(i.)] 
Remark that 
\begin{equation}\label{eq1_proof_Lem2}
   \Big\|\frac{ \partial \widehat L_n(\theta)}{\partial \theta} - \frac{ \partial  L_n(\theta)}{\partial \theta}  \Big\|_\Theta
 \leq
 \sum_{t=1}^{n}\Big\|\frac{\partial \widehat q_{t}(\theta)}{\partial \theta} - \frac{\partial q_{t}(\theta)}{\partial \theta} \Big\|_\Theta.
\end{equation}
Moreover, for all  $\theta \in \Theta$,  
\begin{align}\label{deriv1_q_t}
 \frac{\partial q_{t}(\theta)}{\partial \theta}
 &= 
 -(H^t_\theta)^{-2}\big(2H^t_\theta(Y_t-  f^t_\theta) \frac{\partial f^t_\theta}{\partial\theta}+(Y_t-  f^t_\theta)^2 \frac{\partial H^t_\theta}{\partial\theta}\big) + (H^t_\theta)^{-1}\frac{H^t_\theta}{\partial\theta} \nonumber\\
 &= 
 -2(H^t_\theta)^{-1}(Y_t-  f^t_\theta) \frac{\partial f^t_\theta}{\partial\theta} +(Y_t-  f^t_\theta)^2 \frac{\partial (H^t_\theta)^{-1}}{\partial\theta} + (H^t_\theta)^{-1}\frac{\partial H^t_\theta}{\partial\theta} ,
 \end{align}
 which implies 
  \begin{align*}
\Big| \frac{\partial \widehat q_{t}(\theta)}{\partial \theta} -\frac{\partial  q_{t}(\theta)}{\partial \theta} \Big|
&\leq
2 \Big| (\widehat H^t_\theta)^{-1}(Y_t-  \widehat f^t_\theta) \frac{\partial \widehat f^t_\theta}{\partial\theta}-(H^t_\theta)^{-1}(Y_t-  f^t_\theta) \frac{\partial f^t_\theta}{\partial\theta}\Big| \\
&\hspace{.5cm} +\Big|(Y_t-  \widehat f^t_\theta)^2 \frac{\partial (\widehat H^t_\theta)^{-1}}{\partial\theta} -(Y_t-  f^t_\theta)^2 \frac{\partial (H^t_\theta)^{-1}}{\partial\theta}\Big|
+ \Big|(\widehat H^t_\theta)^{-1}\frac{\partial \widehat H^t_\theta}{\partial\theta} -(H^t_\theta)^{-1}\frac{\partial H^t_\theta}{\partial\theta}\Big| .
\end{align*}
   Using the relation $|a_1b_1c_1 - a_2b_2c_2| \leq |a_1-a_2| |b_2||c_2| + |a_1||b_1-b_2| |c_2| + |a_1||b_1||c_1-c_2| , ~ \forall a_1, a_2, b_1, b_2, c_1, c_2, \in \R$,  we get
 \begin{align}\label{eq2_proof_Lem2}
&\Big\| \frac{\partial \widehat q_{t}(\theta)}{\partial \theta} -\frac{\partial  q_{t}(\theta)}{\partial \theta} \Big\|_\Theta \nonumber\\
&\hspace{.5cm}\leq
2 \Big( \big\|(\widehat H^t_\theta)^{-1}-(H^t_\theta)^{-1}\big\|_\Theta \|Y_t-\widehat f^t_\theta\|_\Theta  \Big\|\frac{\partial f^t_\theta}{\partial\theta}\Big\|_\Theta + \big\|(\widehat H^t_\theta)^{-1}\big\|_\Theta \|\widehat f^t_\theta-f^t_\theta\|_\Theta  \Big\|\frac{\partial \widehat f^t_\theta}{\partial\theta}\Big\|_\Theta \nonumber\\
&\hspace{1cm}
 +\big\|(\widehat H^t_\theta)^{-1}\big\|_\Theta \|Y_t-f^t_\theta\|_\Theta  \Big\|\frac{\partial \widehat f^t_\theta}{\partial\theta}-\frac{\partial f^t_\theta}{\partial\theta}\Big\|_\Theta  \Big) + \|(Y_t-  \widehat f^t_\theta)\|^2_\Theta \Big\|\frac{\partial (\widehat H^t_\theta)^{-1}}{\partial\theta} - \frac{\partial (H^t_\theta)^{-1}}{\partial\theta}\Big\|_\Theta \nonumber\\
 &\hspace{1cm}
 +2|Y_t|\|\widehat f^t_\theta-f^t_\theta\|_\Theta\Big\|\frac{\partial (H^t_\theta)^{-1}}{\partial\theta}\Big\|_\Theta 
 + \big\|(\widehat H^t_\theta)^{-1}\big\|_\Theta \Big\|\frac{\partial\widehat H^t_\theta}{\partial\theta} -\frac{\partial H^t_\theta}{\partial\theta}\Big\|_\Theta +
  \big\|(\widehat H^t_\theta)^{-1}-(H^t_\theta)^{-1}\big\|_\Theta \Big\|\frac{\partial H^t_\theta}{\partial\theta}\Big\|_\Theta \nonumber\\
  &\hspace{.5cm} \leq
2(\underline{h})^{-1} \Big(  \|\widehat f^t_\theta-f^t_\theta\|_\Theta  \Big\|\frac{\partial \widehat f^t_\theta}{\partial\theta}\Big\|_\Theta + \|Y_t-f^t_\theta\|_\Theta  \Big\|\frac{\partial \widehat f^t_\theta}{\partial\theta}-\frac{\partial f^t_\theta}{\partial\theta}\Big\|_\Theta+
 \frac{1}{2} \Big\|\frac{\partial\widehat H^t_\theta}{\partial\theta} -\frac{\partial H^t_\theta}{\partial\theta}\Big\|_\Theta\Big ) \nonumber\\
 &\hspace{1cm} +
 2  \big\|(\widehat H^t_\theta)^{-1}-(H^t_\theta)^{-1}\big\|_\Theta \|Y_t-\widehat f^t_\theta\|_\Theta  \Big\|\frac{\partial f^t_\theta}{\partial\theta}\Big\|_\Theta + \|(Y_t-  \widehat f^t_\theta)\|^2_\Theta \Big\|\frac{\partial (\widehat H^t_\theta)^{-1}}{\partial\theta} - \frac{\partial (H^t_\theta)^{-1}}{\partial\theta}\Big\|_\Theta \nonumber\\
 &\hspace{1cm}
  +2|Y_t|\|\widehat f^t_\theta-f^t_\theta\|_\Theta\Big\|\frac{\partial (H^t_\theta)^{-1}}{\partial\theta}\Big\|_\Theta +
  \big\|(\widehat H^t_\theta)^{-1}-(H^t_\theta)^{-1}\big\|_\Theta \Big\|\frac{\partial H^t_\theta}{\partial\theta}\Big\|_\Theta .
 \end{align}
By applying the H\"{o}lder's inequality to the terms of the right hand side of (\ref{eq2_proof_Lem2}), we have 
 \begin{align*}
&\E \Big[\Big\| \frac{\partial \widehat q_{t}(\theta)}{\partial \theta} -\frac{\partial  q_{t}(\theta)}{\partial \theta} \Big\|_\Theta\Big]\\
&\hspace{.5cm}\leq
C\Bigg[ 
\left(\E[\|\widehat f^t_\theta-f^t_\theta\|^{4}_\Theta]\right)^{1/4}  \Big(\E\Big[\Big\|\frac{\partial \widehat f^t_\theta}{\partial\theta}\Big\|^{4/3}_\Theta\Big]\Big)^{3/4}  + \big(\E[\|Y_t-f^t_\theta\|^{4/3}_\Theta]\big)^{3/4} \Big(\E\Big[\Big\|\frac{\partial \widehat f^t_\theta}{\partial\theta}-\frac{\partial f^t_\theta}{\partial\theta}\Big\|^{4}_\Theta\Big]\Big)^{1/4} \\
  &\hspace{1cm} +
   \Big(\E\Big[\Big\|\frac{\partial\widehat H^t_\theta}{\partial\theta} -\frac{\partial H^t_\theta}{\partial\theta}\Big\|^{2}_\Theta\Big]\Big)^{1/2}
  + \Big(\E\Big[\big\|(\widehat H^t_\theta)^{-1}-(H^t_\theta)^{-1}\big\|^{4}_\Theta\Big]\Big)^{1/4} \left(\E[\|Y_t-\widehat f^t_\theta\|^{4}_\Theta]\right)^{1/4}   
  \Big(\E\Big[\Big\|\frac{\partial f^t_\theta}{\partial\theta}\Big\|^{2}_\Theta\Big]\Big)^{1/2}
  \\
 &\hspace{1.5cm}  +
 \left(\E[\|(Y_t-  \widehat f^t_\theta)\|^{4}_\Theta]\right)^{1/2}   \Big(\E\Big[\Big\|\frac{\partial (\widehat H^t_\theta)^{-1}}{\partial\theta} - \frac{\partial (H^t_\theta)^{-1}}{\partial\theta}\Big\|^{2}_\Theta\Big]\Big)^{1/2}+\left(\E[|Y_t|^{4}]\right)^{1/4} \left(\E[\|\widehat f^t_\theta-f^t_\theta\|^{4}_\Theta]\right)^{1/4} \\
 &\hspace{2cm} 
 \times \Big(\E\Big[\Big\|\frac{\partial (H^t_\theta)^{-1}}{\partial\theta}\Big\|^{2}_\Theta\Big]\Big)^{1/2}
 +
  \left(\E\big[\big\|(\widehat H^t_\theta)^{-1}-(H^t_\theta)^{-1}\big\|^{4}_\Theta \big]\right)^{1/4}  \Big(\E\Big[\Big\|\frac{\partial H^t_\theta}{\partial\theta}\Big\|^{4/3}_\Theta\Big]\Big)^{3/4}
  \Bigg].
   \end{align*}
  Moreover, since $\theta^* \in \Theta(r)$, 
 using \textbf{A}$_i(f_\theta,\Theta)$ and \textbf{A}$_i(M_\theta,\Theta)$ (with $i=0,1$), one can go along similar lines as in Bardet and Wintenberger (2009) to establish the following results: 
 \begin{align}
     & \label{Res1_A1} 
     \bullet ~  \E \Big[ \Big\|\frac{\partial f^{t}_\theta}{\partial\theta}\Big\|^{r}_\Theta +
     \Big\|\frac{\partial \widehat f^{t}_\theta}{\partial\theta}\Big\|^{r}_\Theta 
     + \Big\|\frac{\partial M^{t}_\theta}{\partial\theta}\Big\|^{r}_\Theta + \Big\|\frac{\partial \widehat M^{t}_\theta}{\partial\theta}\Big\|^{r}_\Theta
      + \Big\|\frac{\partial H^{t}_\theta}{\partial\theta}\Big\|^{r/2}_\Theta + \Big\|\frac{ \partial (H^t_\theta)^{-1}}{\partial\theta}\Big\|^{r/2}_\Theta
     \Big]< \infty,\\ 
     &~\nonumber \\
  &  \label{Res2_A1} 
   \bullet ~\left\{
\begin{array}{l}
  \E\Big[\Big\|\frac{\partial \widehat f^t_\theta}{\partial\theta}-\frac{\partial f^t_\theta}{\partial\theta}\Big\|^r_\Theta\Big]
       \leq  C\Big(\sum\limits_{k \geq t}  \big \{\alpha^{(1)}_{k,Y} (f_\theta,\Theta)+\alpha^{(1)}_{k,X} (f_\theta,\Theta)  \big\}\Big)^{r}, 
       \\
   \E\Big[\big\|(\widehat H^t_\theta)^{-1}-(H^t_\theta)^{-1}\big\|^{r}_\Theta\Big]
       \leq  C\Big(\sum\limits_{k \geq t}  \big \{ \alpha^{(0)}_{k,Y} (M_\theta,\Theta)+\alpha^{(0)}_{k,X} (M_\theta,\Theta)  \big \}\Big)^{r} ,
       \\
       \E\Big[\Big\|\frac{\partial \widehat H^t_\theta}{\partial\theta}-\frac{\partial H^t_\theta}{\partial\theta}\Big\|^{r/2}_\Theta\Big]
       \leq  C\Big(\sum\limits_{k \geq t} \big \{ \alpha^{(0)}_{k,Y} (M_\theta,\Theta)+\alpha^{(1)}_{k,X} (M_\theta,\Theta)  \big \}\Big)^{r/2},
        \\
       \E\Big[ \Big\|\frac{\partial (\widehat H^t_\theta)^{-1}}{\partial\theta}-\frac{\partial (H^t_\theta)^{-1}}{\partial\theta}\Big\|^{r/2}_\Theta \Big]
       \leq 
        C\Big(\sum\limits_{k \geq t} \big \{ \alpha^{(0)}_{k,Y} (M_\theta,\Theta)+\alpha^{(1)}_{k,X} (M_\theta,\Theta)  \big \}\Big)^{r/2}.
\end{array}
\right.
\end{align} 
  Thus, using (\ref{Res1_A0}), (\ref{Res2_A0}), (\ref{Res1_A1}) and (\ref{Res2_A1}) with $r=4$, we obtain 
  \begin{align*}
&\E \Big[\Big\| \frac{\partial \widehat q_{t}(\theta)}{\partial \theta} -\frac{\partial  q_{t}(\theta)}{\partial \theta} \Big\|_\Theta\Big]\\
  & \hspace{1cm}\leq
  C\Bigg[ 
\left(\E[\|\widehat f^t_\theta-f^t_\theta\|^{4}_\Theta]\right)^{1/4}    +  \Big(\E\Big[\Big\|\frac{\partial \widehat f^t_\theta}{\partial\theta}-\frac{\partial f^t_\theta}{\partial\theta}\Big\|^{4}_\Theta\Big]\Big)^{1/4} 
   +
   \Big(\E\Big[\Big\|\frac{\partial\widehat H^t_\theta}{\partial\theta} -\frac{\partial H^t_\theta}{\partial\theta}\Big\|^{2}_\Theta\Big]\Big)^{1/2} \\
  &\hspace{2cm}
  + \Big(\E\Big[\big\|(\widehat H^t_\theta)^{-1}-(H^t_\theta)^{-1}\big\|^{4}_\Theta\Big]\Big)^{1/4}   +
  \Big(\E\Big[\Big\|\frac{\partial (\widehat H^t_\theta)^{-1}}{\partial\theta} - \frac{\partial (H^t_\theta)^{-1}}{\partial\theta}\Big\|^{2}_\Theta\Big]\Big)^{1/2}+ \left(\E[\|\widehat f^t_\theta-f^t_\theta\|^{4}_\Theta]\right)^{1/4} \\
 &\hspace{2.5cm} 
 +
  \left(\E\big[\big\|(\widehat H^t_\theta)^{-1}-(H^t_\theta)^{-1}\big\|^{4}_\Theta \big]\right)^{1/4} 
  \Bigg]\\
&\hspace{1cm}\leq
C\sum\limits_{k \geq t}  \Big\{\alpha^{(0)}_{k,Y} (f_\theta,\Theta) +\alpha^{(0)}_{k,X} (f_\theta,\Theta)+\alpha^{(1)}_{k,Y} (f_\theta,\Theta) +\alpha^{(1)}_{k,X} (f_\theta,\Theta) \\
&\hspace{3cm} +\alpha^{(0)}_{k,Y} (M_\theta,\Theta) +\alpha^{(0)}_{k,X} (M_\theta,\Theta) + \alpha^{(1)}_{k,Y} (M_\theta,\Theta)+\alpha^{(1)}_{k,X} (M_\theta,\Theta)\Big\}.
\end{align*}
Therefore, in view of the condition (\ref{cond_th2}), it holds that 
\begin{align*}
\E \Big[\Big\|\frac{\partial \widehat q_{t}(\theta)}{\partial \theta} - \frac{\partial q_{t}(\theta)}{\partial \theta} \Big\|_\Theta\Big] 
 & \leq 
 C\sum\limits_{k \geq t} k^{- \gamma} =C \frac{1}{t^{\gamma-1}}.
 \end{align*}
 By the inequality (\ref{eq1_proof_Lem2}), we deduce
 \begin{align*}
 \E\Big[ \frac{1}{\sqrt{n}} \Big\|\frac{ \partial \widehat L_n(\theta)}{\partial \theta} - \frac{ \partial  L_n(\theta)}{\partial \theta}  \Big\|_\Theta \Big]
 & \leq  
 C\frac{1}{\sqrt{n}}  \sum\limits_{t=1}^{n} \frac{1}{t^{\gamma-1}}= C \frac{1}{\sqrt{n}} (1+n^{2-\gamma}) \limiten 0.
 \end{align*}
 This proves the part (i.) of Lemma \ref{lem2}.\\
 
 \item[(ii.)]  
      This part can be established by using the same arguments as in the proof of Lemma \ref{lem1}.\\
 
 \item[(iii.)] 
 Let us show that $\E\Big[ \Big\|\frac{ \partial^2  q_{t}(\theta)}{\partial \theta_i \partial \theta_j}\Big\|_\Theta\Big]  <\infty$, for all $i,j \in \{1,\ldots, d\}$.\\
 From (\ref{deriv1_q_t}), for any $i,j \in \{1,\ldots, d\}$, we have
 \begin{align*}
 \frac{ \partial^2  q_{t}(\theta)}{\partial \theta \partial_i \theta_j}
 &=
 -2 (H^t_\theta)^{-1} (Y_t-f^t_\theta)\frac{\partial^2  f^t_\theta}{\partial \theta_i \partial \theta_j} + (Y_t-f^t_\theta)^2\frac{\partial^2  (H^t_\theta)^{-1}}{\partial \theta_i \partial \theta_j} 
 -2 (Y_t-f^t_\theta)\Big(\frac{\partial f^t_\theta}{\partial \theta_i} \frac{\partial (H^t_\theta)^{-1}}{\partial \theta_j} +\frac{\partial f^t_\theta}{\partial \theta_j} \frac{\partial (H^t_\theta)^{-1}}{\partial \theta_i}  \Big) \\
 & \hspace{.5cm}
 +2 (H^t_\theta)^{-1}\frac{\partial f^t_\theta}{\partial \theta_i} \frac{\partial f^t_\theta}{\partial \theta_j}
 + \frac{\partial (H^t_\theta)^{-1}}{\partial \theta_j}\frac{\partial H^t_\theta}{\partial \theta_i} +(H^t_\theta)^{-1}\frac{\partial^2 H^t_\theta}{\partial \theta_i \partial \theta_j}.
  \end{align*}
 Therefore, according to (\textbf{A1}), we get 
 \begin{align}\label{eq3_proof_Lem2}
&\Big\|\frac{ \partial^2  q_{t}(\theta)}{\partial \theta \partial_i \theta_j}\Big\|_\Theta \nonumber\\
& \hspace{.2cm} \leq
C
\|(Y_t-f^t_\theta)\|_\Theta \Big(\Big\|\frac{\partial^2  f^t_\theta}{\partial \theta_i \partial \theta_j}\Big\|_\Theta+ \Big\|\frac{\partial f^t_\theta}{\partial \theta_i}\Big\|_\Theta  \Big\|\frac{\partial (H^t_\theta)^{-1}}{\partial \theta_j}\Big\|_\Theta +\Big\|\frac{\partial f^t_\theta}{\partial \theta_j}\Big\|_\Theta  \Big\|\frac{\partial (H^t_\theta)^{-1}}{\partial \theta_i}\Big\|_\Theta  \Big) \nonumber\\
& \hspace{.5cm} +
C\Big( \Big\|\frac{\partial f^t_\theta}{\partial \theta_i}\Big\|_\Theta  \Big\|\frac{\partial f^t_\theta}{\partial \theta_j}\Big\|_\Theta +  \Big\|\frac{\partial^2 H^t_\theta}{\partial \theta_i \partial \theta_j}\Big\|_\Theta\Big) +{\|(Y_t-f^t_\theta)\|}^{2}_\Theta \Big\|\frac{\partial^2  (H^t_\theta)^{-1}}{\partial \theta_i \partial \theta_j}\Big\|_\Theta + \Big\|\frac{\partial (H^t_\theta)^{-1}}{\partial \theta_j}\Big\|_\Theta \Big\|\frac{\partial H^t_\theta}{\partial \theta_i}\Big\|_\Theta.
\end{align}
  Moreover, by  \textbf{A}$_2(f_\theta,\Theta)$ and \textbf{A}$_2(M_\theta,\Theta)$, one can show that 
  $$\E\Big[ \Big\|\frac{\partial^2  f^t_\theta}{\partial \theta_i \partial \theta_j}\Big\|^4_\Theta + \Big\|\frac{\partial^2  H^t_\theta}{\partial \theta_i \partial \theta_j}\Big\|^{2}_\Theta +\Big\|\frac{\partial^2  (H^t_\theta)^{-1}}{\partial \theta_i \partial \theta_j}\Big\|^{2}_\Theta\Big]<\infty.$$
  Thus, by applying the H\"{o}lder's inequality to the terms of the right hand side of (\ref{eq3_proof_Lem2}), it suffices to use (\ref{Res1_A0}) and (\ref{Res1_A1}) to obtain $\E\Big[ \Big\|\frac{ \partial^2  q_{t}(\theta)}{\partial \theta_i \partial \theta_j}\Big\|_\Theta\Big] <\infty$. \\
Since $\E\Big[ \Big\|\frac{ \partial^2  q_{t}(\theta)}{\partial \theta_i \partial \theta_j}\Big\|_\Theta\Big] <\infty$ for all $i,j \in \{1,\ldots, d\}$,  from the  stationarity and ergodicity properties of $\big\{\frac{ \partial^2  q_{t}(\theta)}{\partial \theta \partial \theta'},\ t \in \Z \big\}$ and the uniform strong law of large numbers, it holds that 
\[
\Big\|\frac{1}{n} \sum_{t=1}^{n}\frac{ \partial^2 q_{t}(\theta)}{\partial \theta \partial \theta'} - \E\Big( \frac{ \partial^2  q_{0}(\theta)}{\partial \theta \partial \theta'}\Big)  \Big\|_\Theta \limitepsn 0 .
\]
This completes the proof of Lemma \ref{lem2}. 
$~~~~~~~~~~~~~~~~~~~~~~~~~~~~~~~~~~~~~~~~~~~~~~~~~~~~~~~~~~~~~~~~~~~~~~~~~~~~~~~~~~~~~~~~~~~~~~~~
\blacksquare$ 
\end{enumerate} 
 \medskip
    
\noindent
The following lemma is also needed.
\begin{lem}\label{lem3} 
 Assume that the conditions of Theorem \ref{th2} hold. Then
 \rm
\begin{enumerate}
    \item [(i.)] $\big\{\frac{ \partial  q_{t}(\theta^*)}{\partial \theta}|\mathcal{F}_{t-1},\, t \in \Z \big\}$ is a stationary ergodic martingale difference sequence with covariance matrix $G$,  
	\item [(ii.)]  $ -\frac{1}{n}\frac{\partial^2}{\partial \theta \partial \theta'} \widehat L_n( \tilde  \theta_{n})  \limitepsn F$, for any sequence $(\tilde  \theta_{n})_{n\geq1}$ with values in $\Theta$  and satisfying $\tilde  \theta_{n} \limitepsn \theta^*$,
\end{enumerate}
where  $G$ and $F$ are defined in (\ref{def_matrix_F_G}). 
\end{lem}

\medskip

\noindent
\emph{\bf Proof of Lemma \ref{lem3}}
\begin{enumerate}
 \item [(i.)] Recall that $G= \E\Big[\frac{\partial q_{t}(\theta^*)}{\partial \theta} \frac{\partial q_{t}(\theta^*)}{\partial \theta'}\Big]$ and 
 that for all $\theta \in \Theta$,
 \[
 \frac{\partial q_{t}(\theta)}{\partial \theta}=
 -2(H^t_\theta)^{-1}(Y_t-  f^t_\theta) \frac{\partial f^t_\theta}{\partial\theta} -\Big(\frac{Y_t-  f^t_\theta}{H^t_\theta}\Big)^2  \frac{\partial H^t_\theta}{\partial\theta} + (H^t_\theta)^{-1}\frac{\partial H^t_\theta}{\partial\theta}.
 \]
 Since the functions $f^t_\theta$, $H^t_\theta$, $\frac{\partial f^t_\theta}{\partial \theta}$ and $\frac{\partial H^t_\theta}{\partial \theta}$ are $\mathcal{F}_{t-1}$-measurable, we have
 \[
  \E\Big[\frac{ \partial  q_{t}(\theta^*)}{\partial \theta} \big| \mathcal{F}_{t-1}\Big] =
  -(H^t_{\theta^*})^{-1}\frac{\partial H^t_{\theta^*}}{\partial\theta} \left((H^t_{\theta^*})^{-1}\E \left[(Y_t-  f^t_{\theta^*})^2 |\mathcal{F}_{t-1} \right]   -1\right)=0, 
  \]
  which shows that (i.) holds.
  
  \item [(ii.)]  Let $(\tilde  \theta_{n} )_{n \in \N}$ be a sequence satisfying $\tilde  \theta_{n} \limitepsn \theta^*$. 
  For any $i,j =1,\ldots,d$, we have
		 \begin{align*}
		& \Big| \frac{1}{n} \sum\limits_{t =1}^{n}\frac{\partial}{\partial \theta_j \partial \theta_i} q_{t}(\tilde \theta_{n}) 
		 -\E\Big(\frac{\partial}{\partial \theta_j \partial \theta_i} q_{0}( \theta^*)\Big)\Big| \\
		 &\leq
		 \Big| \frac{1}{n} \sum\limits_{t =1}^{n}\frac{\partial}{\partial \theta_j \partial \theta_i} q_{t}(\tilde \theta_{n}) 
		 -\E\Big(\frac{\partial}{\partial \theta_j \partial \theta_i} q_{0}(\tilde \theta_{n})\Big)\Big|
	 +
		 \Big|\E\Big(\frac{\partial}{\partial \theta_j \partial \theta_i} q_{0}( \tilde \theta_{n})\Big)
		 -\E\Big(\frac{\partial}{\partial \theta_j \partial \theta_i} q_{0}(\theta^*)\Big)\Big|\\
		 &\leq
		 \Big\| \frac{1}{n} \sum\limits_{t =1}^{n}\frac{\partial}{\partial \theta_j \partial \theta_i} q_{t}(\theta) 
		 -\E\Big(\frac{\partial}{\partial \theta_j \partial \theta_i} q_{0}( \theta)\Big) \Big\|_\Theta
		 +
		 \Big|\E\Big(\frac{\partial}{\partial \theta_j \partial \theta_i} q_{0}( \tilde \theta_{n})\Big)
		 -\E\Big(\frac{\partial}{\partial \theta_j \partial \theta_i} q_{0}(\theta^*)\Big)\Big|\\
		  &\limiten 0~~\text{(by virtue of Lemma \ref{lem2} (iii.))}.
		 \end{align*}
%
 Thus,
\[
 -\frac{1}{n}\frac{\partial^2}{\partial \theta \partial \theta'}L_{n}(\tilde \theta_{n})  
= \frac{1}{n} \sum\limits_{t =1}^{n}\frac{\partial^2}{\partial \theta \partial \theta'}q_{t}(\tilde \theta_{n}) 
 \limitepsn \E\Big(\frac{\partial^2}{\partial \theta \partial \theta'} q_{0}( \theta^*)\Big)=F.
\]
We conclude the proof of the part (ii.) by using  Lemma \ref{lem2} (ii.).
$~~~~~~~~~~~~~~~~~~~~~~~~~~~~~~~~~~~~~~ 
\blacksquare$ 
 \end{enumerate}

\medskip

\noindent 
Now, we use the results of Lemma \ref{lem2} and \ref{lem3} to prove  the first part of Theorem \ref{th2}.
 The second part can be established by going along similar lines as in Kengne (2021).\\ 
By applying a second-order Taylor expansion to the function $\theta \mapsto \widehat L_n(\theta)$, for all $\theta \in \Theta$, there exists $\tilde \theta$ between $\theta$ and $\theta^*$ such that
\begin{equation}\label{eq1_proof_th2}
 \frac{1}{n} \left\{\widehat L_n(\theta) - \widehat L_n(\theta^*)\right\}
=
\frac{1}{n}\frac{\partial  L_n( \theta^*)}{\partial \theta'}(\theta-\theta^*) -\frac{1}{2}(\theta-\theta^*)' F (\theta-\theta^*) +
R_n(\theta),
\end{equation}
 where
\[
R_n(\theta)= \frac{1}{n} \Big\{\frac{\partial \widehat L_n( \theta^*)}{\partial \theta'}  -\frac{\partial L_n( \theta^*)}{\partial \theta'}\Big\}
(\theta-\theta^*)
+
\frac{1}{2}(\theta-\theta^*)' \Big(\frac{1}{n}\frac{\partial^2}{\partial \theta \partial \theta'} \widehat L_n( \tilde  \theta)+F \Big) (\theta-\theta^*) .
\]
Let us define the vector
\[
Z_n= F^{-1} \frac{1}{\sqrt{n}}\frac{\partial L_n( \theta^*)}{\partial \theta}. 
\]
Then, we can rewrite (\ref{eq1_proof_th2}) as  
\begin{equation}\label{eq2_proof_th2}
\frac{1}{n} \left\{\widehat L_n(\theta) - \widehat L_n(\theta^*)\right\} = \frac{1}{2n} \left\|Z_n \right\|^2_F-\frac{1}{2n} \left\|Z_n-\sqrt{n}(\theta-\theta^*) \right\|^2_F +R_n(\theta).
\end{equation}
Define also 
 \[
 \theta_{Z_n}= \text{arg} \inf_{\theta \in \Theta}\left\|Z_n-\sqrt{n}(\theta-\theta^*) \right\|_F.
 \]
Then, by (\ref{eq0_normal}), for $n$ large enough, we have
\begin{equation*}
\sqrt{n}(\theta_{Z_n}-\theta^*)=Z^{\mathcal C}_n,
\end{equation*}
where $Z^{\mathcal C}_n$ is  the $F$-projection of $Z_n$ on $\mathcal C$. 
Using this relation and the definition of $\theta_{Z_n}$,  we have 
\[
\big\|Z_n-\sqrt{n}(\widehat \theta_{n}-\theta^*) \big\|^2_F -\big\|Z_n - Z^{\mathcal C}_n\big\|^2_F=
\big\|Z_n-\sqrt{n}(\widehat \theta_{n}-\theta^*) \big\|^2_F -\left\|Z_n-\sqrt{n}(\theta_{Z_n}-\theta^*) \right\|^2_F \geq 0.
\]
Furthermore, from (\ref{eq2_proof_th2}) and the definition of $\widehat \theta_{n}$,  it holds that
\begin{align*}
\big\|Z_n-\sqrt{n}(\widehat \theta_{n}-\theta^*) \big\|^2_F -\left\|Z_n-\sqrt{n}(\theta_{Z_n}-\theta^*) \right\|^2_F
&=
\{\widehat L_n(\theta_{Z_n}) - \widehat L_n(\widehat \theta_{n})\} + 2n\{R_n(\widehat \theta_{n}) - R_n(\theta_{Z_n})\}\\
&\leq 
2n\{R_n(\widehat \theta_{n}) - R_n(\theta_{Z_n})\}.
\end{align*}
Therefore, 
\begin{equation}\label{eq3_proof_th2}
\left|\big\|Z_n-\sqrt{n}(\widehat \theta_{n}-\theta^*) \big\|^2_F -\big\|Z_n - Z^{\mathcal C}_n\big\|^2_F  \right|
\leq
2n\{R_n(\widehat \theta_{n}) - R_n(\theta_{Z_n})\}.
\end{equation}
Let us consider the following Lemma.
\begin{lem}\label{lem4} 
 Assume that the conditions of Theorem \ref{th2} hold. Then
 \[
 n \{R_n(\widehat \theta_{n})  - R_n(\theta_{Z_n})\}  =o_P(1).
 \]
\end{lem}
By Lemma \ref{lem4} and (\ref{eq3_proof_th2}), it follows that 
\begin{equation}\label{eq4_proof_th2}
\big\|Z_n-\sqrt{n}(\widehat \theta_{n}-\theta^*) \big\|^2_F -\big\|Z^{\mathcal C}_n - Z_n\big\|^2_F  
=o_P(1).
\end{equation}
Moreover, according to the equivalent definition of the $F$-orthogonal projection in (\ref{def_project_orth}), we get   
\begin{align*}
\big\|Z_n-\sqrt{n}(\widehat \theta_{n}-\theta^*) \big\|^2_F  
& =
\big\|Z^{\mathcal C}_n-\sqrt{n}(\widehat \theta_{n}-\theta^*) \big\|^2_F +
\big\|Z^{\mathcal C}_n-Z_n\big\|^2_F 
- 2\left\langle Z^{\mathcal C}_n-\sqrt{n}(\widehat \theta_{n}-\theta^*), Z^{\mathcal C}_n-Z_n  \right\rangle\\
&\geq
\big\|Z^{\mathcal C}_n-\sqrt{n}(\widehat \theta_{n}-\theta^*) \big\|^2_F +
\big\|Z^{\mathcal C}_n-Z_n\big\|^2_F .
\end{align*}
Therefore, from (\ref{eq4_proof_th2}), we obtain
\begin{equation}\label{oP_Zn_C}
\big\|Z^{\mathcal C}_n-\sqrt{n}(\widehat \theta_{n}-\theta^*) \big\|^2_F  
\leq
\big\|Z_n-\sqrt{n}(\widehat \theta_{n}-\theta^*) \big\|^2_F -\big\|Z^{\mathcal C}_n - Z_n\big\|^2_F  
=o_P(1).
\end{equation}
 Now, using Lemma \ref{lem3} (i.), we apply the central limit theorem for the stationary ergodic martingale difference sequence $\big\{\frac{ \partial  q_{t}(\theta^*)}{\partial \theta}|\mathcal{F}_{t-1},\, t \in \Z \big\}$.  
It follows that
\begin{equation}\label{Normal_d_L}
\frac{1}{\sqrt{n}}\frac{\partial L_{n}(\theta^*)}{\partial \theta} =\frac{1}{\sqrt{n}}\sum\limits_{t =1}^{n}\frac{\partial q_{t}(\theta^*)}{\partial \theta}
\limiteloin \mathcal{N}_d \left(0,G \right),
\end{equation}
and thus
\begin{align}\label{eq5_proof_th2}
   Z_n=F^{-1} \frac{1}{\sqrt{n}}\frac{\partial L_n( \theta^*)}{\partial \theta}
     \limiteloin Z \sim 
      \mathcal{N}_d \left(0,F^{-1} G F^{-1} \right).
    \end{align}
Hence, $Z^{\mathcal C}_n \limiteloin Z^{\mathcal C}$. From this, it suffices to use (\ref{oP_Zn_C}) to conclude the proof of Theorem \ref{th2}.
     \begin{flushright}
$\blacksquare$
\end{flushright}

\noindent 
{\bf Proof of Lemma \ref{lem4}.}\\
Recall that 
\[
R_n(\theta)= \frac{1}{n} \Big\{\frac{\partial \widehat L_n( \theta^*)}{\partial \theta'}  -\frac{\partial  L_n( \theta^*)}{\partial \theta'}\Big\}
(\theta-\theta^*)
+
\frac{1}{2}(\theta-\theta^*)' \Big(\frac{1}{n}\frac{\partial^2}{\partial \theta \partial \theta'} \widehat L_n( \tilde  \theta)+F \Big) (\theta-\theta^*) .
\]
According to Lemmas \ref{lem2} (i.) and \ref{lem3} (ii.), when $\tilde \theta_n -\theta^*=o_P(1)$, we have 
\begin{equation}\label{eq1_proof_lem4}
nR_n(\tilde \theta_n)=o_P(\sqrt{n}(\tilde \theta_{n}-\theta^*))+ o_P(n  \|\tilde \theta_{n}-\theta^*\|^2).
\end{equation}
This implies 
\begin{equation}\label{eq2_proof_lem4}
nR_n(\tilde \theta_n)=o_P(1) ~\text{ when }~ \sqrt{n}(\tilde \theta_{n}-\theta^*)=O_P(1).
\end{equation}
It comes from the definition of $\theta_{Z_n}$ that  
\[
 \left\| \sqrt{n}(\theta_{Z_n}-\theta^*)\right\|_F \leq \left\| \sqrt{n}(\theta_{Z_n}-\theta^*) -Z_n\right\|_F +\left\| Z_n\right\|_F  \leq 2\left\| Z_n\right\|_F .
\]
Moreover, the convergence in (\ref{eq5_proof_th2}) implies $\left\| Z_n\right\|_F=O_P(1)$; and consequently,   $\sqrt{n}(\theta_{Z_n}-\theta^*)=O_P(1)$. Thus,   $nR_n(\theta_{Z_n})=o_P(1)$ by  virtue (\ref{eq2_proof_lem4}).
\medskip
   
\noindent
    We now show that, it also holds  $nR_n(\widehat \theta_{n})=o_P(1)$. From (\ref{eq2_proof_th2}), we have
\begin{equation*}
\left\|Z_n \right\|^2_F-\big\|Z_n-\sqrt{n}(\widehat \theta_{n}-\theta^*) \big\|^2_F +2nR_n(\widehat \theta_{n}) =2\{\widehat L_n(\widehat \theta_{n}) - \widehat L_n(\theta^*)\} \geq 0,
\end{equation*}
where the inequality holds since $\widehat{\theta}_{n}= \underset{\theta\in \Theta}{\text{argmax}} (\widehat L_n(\theta) )$. 
Thus, it holds that 
\begin{align*}
\big \|\sqrt{n}(\widehat{\theta}_{n}-\theta^*) \big\|^2_F
& \leq
2 \big(\big\|Z_n-\sqrt{n}(\widehat{\theta}_{n}-\theta^*)\big\|^2_F +\left\|Z_n\right\|^2_F\big)\\
& \leq 4\left\|Z_n\right\|^2_F +4 nR_n(\widehat \theta_{n}). 
\end{align*}
Furthermore, since $\widehat{\theta}_{n} \limitepsn \theta^*$,  by (\ref{eq1_proof_lem4}), it follows that $nR_n(\widehat \theta_{n}) =o_P\big( \big\|\sqrt{n}(\widehat{\theta}_{n}-\theta^*)\big\|^2_F\big)$. Consequently,  $\sqrt{n}(\widehat{\theta}_{n}-\theta^*)=O_P(1)$, and $nR_n(\widehat \theta_{n})=o_P(1)$ holds according to  (\ref{eq2_proof_lem4}). 
 This achieves the proof of the lemma.
 \begin{flushright}
 $~~~~~~\blacksquare$ 
\end{flushright}

 \subsection{Proof of Theorem \ref{th3}}
Under $H_0$, we have  $\Gamma \widehat\theta_{n}-\vartheta_0= \Gamma ( \widehat\theta_{n}-\theta^*)$. 
Then, we get
\begin{align}\label{eq_proofth3}
W_n &= n(\Gamma \widehat\theta_{n}-\vartheta_0)' (\Gamma  \widehat \Sigma_n \Gamma ')^{-1}(\Gamma \widehat\theta_{n}-\vartheta_0)\nonumber\\
&=
n ( \widehat\theta_{n}-\theta^*)' \Gamma'(\Gamma  \widehat \Sigma_n \Gamma ')^{-1} \Gamma ( \widehat\theta_{n}-\theta^*) \nonumber\\
& =
\sqrt{n} ( \widehat\theta_{n}-\theta^*)' \Gamma'(\Gamma \Sigma \Gamma ')^{-1} \Gamma \sqrt{n} ( \widehat\theta_{n}-\theta^*)\nonumber \\
& \hspace{5cm}+
\sqrt{n} ( \widehat\theta_{n}-\theta^*)' \Gamma'\left((\Gamma \widehat \Sigma_n  \Gamma ')^{-1} -(\Gamma \Sigma \Gamma ')^{-1}\right) \Gamma \sqrt{n} ( \widehat\theta_{n}-\theta^*).
\end{align}
Recall that, by Theorem \ref{th3}, we have $\sqrt{n} ( \widehat\theta_{n}-\theta^*) \limiteloin  Z^{\mathcal C}~$ with $~Z \sim \mathcal{N}_d \left(0,\Sigma \right)$.\\
 Furthermore, $(\Gamma \widehat \Sigma_n  \Gamma ')^{-1} -(\Gamma \Sigma \Gamma ')^{-1} =o_P(1)$.  
Thus, from  (\ref{eq_proofth3}), it holds that
\begin{align*}
W_n 
&=
\sqrt{n} ( \widehat\theta_{n}-\theta^*)' \Gamma'(\Gamma \Sigma \Gamma ')^{-1} \Gamma \sqrt{n} ( \widehat\theta_{n}-\theta^*) + o_P(1)\\
&~~\limiteloin  
 (\Gamma Z^{\mathcal C})'(\Gamma \Sigma \Gamma')^{-1}\Gamma Z,
\end{align*}
which  establishes the theorem. 
 \begin{flushright}
$\blacksquare$
\end{flushright}

\noindent 
{\bf Proof of Corollary \ref{corol1}.}\\
When $\theta^* \in \overset{\circ}{\Theta}$, we have  $W_n \limiteloin  Z' \Gamma'(\Gamma \Sigma \Gamma')^{-1}\Gamma Z={\left\| U\right\|}^2$ 
with $U=(\Gamma \Sigma \Gamma')^{-1/2}\Gamma Z$ and $Z \sim \mathcal{N}_d \left(0,\Sigma \right)$. 
Since $\Sigma$ is symmetric, the vector $U$  follows a multivariate Gaussian distribution with mean $0$ and covariance matrix $I_{d_0}$,  where $I_{d_0}$ is the identity matrix of size $d_0$.
Therefore, all components of $U$ are independent, standard normal distributed random variables. This leads to the conclusion.
\begin{flushright}
$\blacksquare$
\end{flushright}

\subsection{Proof of Theorem \ref{th4}}
Consider the following lemma.
\begin{lem}\label{lem5} 
 Assume that the conditions of Theorem \ref{th4} hold. Then
 \[
 \frac{1}{\sqrt{n \log \log n}} \Big\|\frac{ \partial \widehat L_n(\theta)}{\partial \theta} - \frac{ \partial  L_n(\theta)}{\partial \theta}  \Big\|_\Theta   \limitepsn 0.
 \]
  \end{lem}

\noindent
\emph{\bf Proof of Lemma \ref{lem5}.}\\
Using the inequality (\ref{eq1_proof_Lem2}) and Corollary 1 of of Kounias and Weng (1969), it suffices to show that 
\begin{equation}\label{Kounias_lem5}
 \sum_{k \geq 2}^{}\frac{1}{\sqrt{k \log \log k}}\E \left[\Big\|\frac{\partial \widehat q_{k}(\theta)}{\partial \theta} - \frac{\partial q_{k}(\theta)}{\partial \theta} \Big\|_\Theta\right] < \infty .
\end{equation}
In the proof of Lemma \ref{lem3}, we have established that
\begin{align*}
\Big\|\frac{\partial \widehat q_{k}(\theta)}{\partial \theta} - \frac{\partial q_{k}(\theta)}{\partial \theta} \Big\|_\Theta
&\leq
C\sum\limits_{j \geq k}  \Big\{\alpha^{(0)}_{j,Y} (f_\theta,\Theta) +\alpha^{(0)}_{j,X} (f_\theta,\Theta)+\alpha^{(1)}_{j,Y} (f_\theta,\Theta) +\alpha^{(1)}_{k,X} (f_\theta,\Theta) \\
 & \hspace{3cm}
 +\alpha^{(0)}_{j,Y} (M_\theta,\Theta) +\alpha^{(0)}_{j,X} (M_\theta,\Theta) + \alpha^{(1)}_{j,Y} (M_\theta,\Theta)+\alpha^{(1)}_{j,X} (M_\theta,\Theta)\Big\}\\
 &=
 C \sum\limits_{j \geq k}  \sum_{i = 0}^{1}
  \left\{ \alpha^{(i)}_{j,Y}(f_\theta,\Theta)+\alpha^{(i)}_{j,X}(f_\theta,\Theta)+\alpha^{(i)}_{j,Y}(M_\theta,\Theta)+\alpha^{(i)}_{j,X}(M_\theta,\Theta)\right\}.
\end{align*} 
Then, from the condition  (\ref{cond_th4}), we obtain
\begin{multline*}
\sum_{k \geq 2}^{}\frac{1}{\sqrt{k \log \log k}}
\E \left[\Big\|\frac{\partial \widehat q_{k}(\theta)}{\partial \theta} - \frac{\partial q_{k}(\theta)}{\partial \theta} \Big\|_\Theta\right]\\
\leq
  \sum_{k \geq 2} \frac{1}{\sqrt{k \log \log k}}\sum_{j \geq k} \sum_{i = 0}^{1}
  \big\{ \alpha^{(i)}_{j,Y}(f_\theta,\Theta)+\alpha^{(i)}_{j,X}(f_\theta,\Theta)+\alpha^{(i)}_{j,Y}(M_\theta,\Theta)+\alpha^{(i)}_{j,X}(M_\theta,\Theta)
   \big\} <
   \infty.
\end{multline*}
Hence, (\ref{Kounias_lem5}) is satisfied, and Lemma \ref{lem5} holds.
 \begin{flushright}
$\blacksquare$
\end{flushright} 

\medskip

\noindent
 Let us prove the part (i.) of the theorem.
\begin{enumerate}
	\item [(i.)]
 We have 
\begin{equation*}
 P(\widehat m_n=m^*)=1-P(\widehat m_n \supsetneq m^*) -P(\widehat m_n \nsupseteq m^*).
\end{equation*}
Therefore, it suffices to show that 
\begin{equation}\label{cond_proof_th4}
\lim\limits_{n \rightarrow \infty} P (\widehat m_n \supsetneq  m^*) = \lim\limits_{n \rightarrow \infty}P (\widehat m_n \nsupseteq m^*)=0.
\end{equation} 

1. Let $m\in {\cal M}$ such as $m  \supsetneq  m^*$. 
We have, 
\begin{equation}\label{C_m_star_C_m_log_log_n}
 \frac{1}{ \sqrt{ \log \log n }} \big( \widehat{C}(m^*) - \widehat{C}(m) \big) = \frac{2}{ \sqrt{\log \log n }} \big(\widehat{L}_n\big(\widehat{\theta}(m)\big)-\widehat{L}_n\big(\widehat{\theta}(m^*)\big) - \frac{\kappa_n}{ \sqrt{\log \log n }}  ( |m|-|m^*|). 
\end{equation}
Let us establish that
\begin{equation}\label{L_m_star_L_m_log_log_n}
 \frac{1}{ \sqrt{\log \log n } } \big(\widehat{L}_n\big(\widehat{\theta}(m)\big)-\widehat{L}_n\big(\widehat{\theta}(m^*)\big)\big) = O_P(1) .
\end{equation}
From the Taylor expansion of $\widehat{L}_n$, we can find $\overline{\theta}(m)$ between $\widehat{\theta}(m)$ and $\theta^*$ such that
%
\begin{equation}\label{Taylor_widehat_L2}
\widehat{L}_n\big(\widehat{\theta}(m)\big)-\widehat{L}_n\big(\theta^*\big) = \dfrac{\partial L_n( \theta^* )}{\partial \theta}\big( \widehat{\theta}(m) - \theta^* \big) -
 \frac{1}{2}\sqrt{n}\big( \widehat{\theta}(m) - \theta^* \big)' F(\theta^*,m)\sqrt{n} \big( \widehat{\theta}(m) - \theta^* \big) + nR'_n(m),
\end{equation}
 where
\begin{multline*}
R'_n(m) = 
\frac{1}{n} \Big\{\frac{\partial \widehat L_n( \theta^*)}{\partial \theta'}  -\frac{\partial L_n( \theta^*)}{\partial \theta'}\Big\}
\big(\widehat \theta(m)-\theta^* \big)\\
+
\frac{1}{2}\big(\widehat \theta(m)-\theta^*\big)' \Big(\frac{1}{n}\frac{\partial^2}{\partial \theta \partial \theta'} \widehat L_n(  \overline  \theta(m))+F(\theta^*,m) \Big) \big(\widehat \theta(m)-\theta^*\big)
\end{multline*}
and
\begin{flalign*}
&  F(\theta^*,m)= \Big(\E \Big[ \dfrac{ \partial^2 q_0(\theta^*)}{ \partial\theta_i \partial\theta_j}  \Big] \Big)_{i,j \in m}. & &
\end{flalign*}
Moreover, since $\widehat{\theta}(m),~ \overline{\theta}(m) \limitepsn \theta^*$, in this case of overfitting, the same arguments as in the proof of Lemma   \ref{lem3} (ii.) lead to
\begin{equation*}
 -\dfrac{1}{n} \dfrac{\partial^2 \widehat{L}_n\big( \overline{\theta}(m) \big)}{\partial \theta \partial \theta'} \limitepsn F(\theta^*,m).
\end{equation*}     
%
%
%
%
%
Then, one can show as in the proof of Theorem \ref{th2} that $nR'_n(m) = o_P(1)$. Also, we have $\sqrt{n} \big( \widehat{\theta}(m) - \theta^* \big)= O_P(1)$. 
 In addition, $\big\{\frac{ \partial  q_{t}(\theta^*)}{\partial \theta}|\mathcal{F}_{t-1},\, t \in \Z \big\}$ is a stationary ergodic square integrable  martingale difference sequence (see above). Hence,  from the law of iterative logarithm for martingales (see for instance \cite{Stout1970, Stout1974}), we get,
\[  \frac{1}{ \sqrt{n \log \log n }} \dfrac{\partial L_n( \theta^* )}{\partial \theta} = O(1). \]
Thus, we have from (\ref{Taylor_widehat_L2}), 
\begin{align}\label{Taylor_widehat_L3}
 &\nonumber \frac{1}{ \sqrt{ \log \log n } } \big(\widehat{L}_n\big(\widehat{\theta}(m)\big)-\widehat{L}_n\big(\theta^* \big)	\big) \\ 
  &=   \frac{1}{ \sqrt{n \log \log n }} \dfrac{\partial L_n( \theta^* )}{\partial \theta} \sqrt{n}\big( \widehat{\theta}(m) - \theta^* \big) + \frac{1}{2\sqrt{\log \log n }} \sqrt{n} \big( \widehat{\theta}(m) - \theta^* \big)' F(\theta^*,m) \sqrt{n} \big( \widehat{\theta}(m) - \theta^* \big) \nonumber  \\
 & \hspace{12.01cm} + \dfrac{1}{ \sqrt{ \log \log n } } nR'_n(m)\nonumber \\
 &  = O(1) O_P(1) + o(1) O_P(1)  O_P(1) + o_P(1) = O_P(1).
 \end{align} 
 By using the same arguments with $m=m^*$, we get 
 \begin{equation}\label{Taylor_widehat_L4}
   \frac{1}{ \sqrt{ \log \log n } } \big(\widehat{L}_n\big(\widehat{\theta}(m^*)\big)- \widehat{L}_n\big(\theta^*) \big) = O_P(1).
 \end{equation} 
 Hence, (\ref{L_m_star_L_m_log_log_n}) holds from (\ref{Taylor_widehat_L3})  and (\ref{Taylor_widehat_L4}). 

 \medskip

\noindent
 Therefore, since $\kappa_n / \sqrt{\log \log n} \limiten \infty$ and $|m| > |m^*|$, then (\ref{C_m_star_C_m_log_log_n}) and (\ref{L_m_star_L_m_log_log_n}) lead to
\begin{equation*}\label{C_m_star_C_m_log_log_n_conv_infty}
  \frac{1}{ \sqrt{\log \log n }} \big( \widehat{C}(m^*) - \widehat{C}(m) \big) \limiteproban -\infty ~ ~  . 
\end{equation*} 
This implies that, for large $n$,
\begin{equation*}\label{C_m_star_C_m_n_large}
   \widehat{C}(m) - \widehat{C}(m^*) > 0 
\end{equation*} 
with probability one; that is, $P(\widehat{m}_n  \supsetneq  m^*) \limiten 0$.

\medskip

2. Let $m\in {\cal M}$ such as $m  \nsupseteq  m^*$. 
  We have,
  \begin{equation}\label{C_m_star_C_m_n}
 \frac{1}{  n} \big( \widehat{C}(m^*) - \widehat{C}(m) \big) = \frac{2}{ n} \big(\widehat{L}_n\big(\widehat{\theta}(m)\big)-\widehat{L}_n\big(\widehat{\theta}(m^*)\big)\big) - \frac{\kappa_n}{ n}  ( |m|-|m^*|). 
\end{equation}
Using the same arguments in the proof of Theorem 3.1 of Bardet \textit{et al.} (2020), we get
\begin{equation*} 
  \frac{1}{ n} \big(\widehat{L}_n\big(\widehat{\theta}(m)\big)-\widehat{L}_n\big(\widehat{\theta}(m^*)\big) \big) 
  = L(\theta^*(m)) - L(\theta^*) + o(1) ~ ~ a.s., 
 \end{equation*}
where $L(\theta) = - \E[q_0(\theta)]$, for all $\theta \in \Theta$.
Note that, the function $L: \Theta \rightarrow \R$ has a unique maximum at $\theta^*$ (see the proof of Theorem \ref{th1}).
Since $m  \nsupseteq  m^*$, it holds that $\theta^* \notin \Theta(m)$; and consequently,
 $L(\theta^*(m)) - L(\theta^*) < 0 ~ ~ a.s.$.
Thus, according to (\ref{C_m_star_C_m_n}) and since $\kappa_n /n \limiten 0$, we get
\begin{equation*}
 \lim_{n \rightarrow \infty} \frac{1}{ n} \big( \widehat{C}(m^*) - \widehat{C}(m) \big) < 0 ~ ~ ~ a.s. ~ \text{ and } ~ ~  \widehat{C}(m)-\widehat{C}(m^*) >0  ~ ~ ~ a.s. ~ \text{ for large } n.
\end{equation*} 
 This implies that $P(\widehat{m}_n  \nsupseteq  m^*) \limiten 0$.
 Hence, the condition (\ref{cond_proof_th4}) holds; and the part (i.) of the theorem  is established.
 
  \item [(ii.)] 
  Let $m\in {\cal M}$ such as $m  \supsetneq  m^*$. We have 
  \begin{equation*}
   \frac{1}{ \log \log n } \big( \widehat{C}(m) - \widehat{C}(m^*) \big) = \frac{2}{ \log \log n } \big(\widehat{L}_n\big(\widehat{\theta}(m^*)\big)-\widehat{L}_n\big(\widehat{\theta}(m)\big)\big) + \frac{\kappa_n}{ \log \log n}  ( |m|-|m^*|). 
  \end{equation*}

Moreover, from the same arguments as in the proof of Theorem 3.1 in Kengne (2021), one can show that 
\[
\frac{1}{ \log \log n  } \big(\widehat{L}_n\big(\widehat{\theta}(m^*)\big)-\widehat{L}_n\big(\widehat{\theta}(m)\big)\big) =  O(1)~~a.s.
\]
Thus, 
 we can find a constant $c$ such that if $\underset{n\rightarrow \infty}{\liminf} \kappa_n/\log \log n>c$, then
 \begin{equation*}
\liminf_{n\rightarrow \infty} \frac{1}{ \log \log n} \big( \widehat{C}(m) - \widehat{C}(m^*) \big) >0~~a.s.
 \end{equation*}
This implies that
 \begin{equation}\label{eq1_proof_ii_th4}
\widehat{C}(m) - \widehat{C}(m^*)>0~~a.s.~~ \text{for large } n.
 \end{equation}
 
 Note that, the inequality (\ref{eq1_proof_ii_th4}) also holds when $m  \nsupseteq  m^*$ (see the part 2. of the proof of (i.)).
 Hence, we deduce that 
$ \widehat{m}_n= \underset{m \in \mathcal{M}}{\text{argmin}} ~ \widehat{C}(m) = \underset{m \in \mathcal{M}}{\text{argmin}} \big( \widehat{C}(m) - \widehat{C}(m^*) \big) \limitepsn m^*$; which establishes the strong consistency of $ \widehat{m}_n$.

 \item [(iii.)] Using Lemma \ref{lem5}, this part can be proved by going along similar lines as in Kengne (2021).
 
 \end{enumerate}
 \begin{flushright}
$\blacksquare$ 
\end{flushright}

 \section*{Acknowledgements}
  The authors are very grateful to the Editor, the Associate Editor and the anonymous Referee for many relevant suggestions and comments which helped to improve the contents of this article.

 \end{document}